\documentclass[a4paper,english,12pt]{article}

\usepackage{amssymb}
\usepackage{epsf}
\usepackage{graphicx}
\usepackage{hyperref}
\usepackage{amssymb}
\usepackage{amsmath}
\usepackage{amsfonts}
\usepackage{dsfont}
\usepackage[T1]{fontenc}
\usepackage{wasysym}
\usepackage{amssymb}
\usepackage{stmaryrd}
\usepackage{mathrsfs}
\usepackage{pbox} % pbox ?
\usepackage{booktabs}
\usepackage{authblk}

\def\RH{2006RH$_{120}$ }
\def\R{\mbox{I\hspace{-.15em}R}}

\title{Rendezvous Missions to Temporarily-Captured Near Earth Asteroids}
\author[1]{S. Brelsford}
\author[1]{M. Chyba}
\author[2]{T.Haberkorn}
\author[1]{G.Patterson}
\affil[1]{Department of Mathematics, University of Hawaii, Honolulu, Hawaii, 96822}
\affil[2]{MAPMO-F\'ed\'eration Denis Poisson, University of Orl\'eans, 45067 Orl\'eans, France}

\begin{document}

\maketitle
% \footnotetext{Department of Mathematics, University of Hawaii, Honolulu, Hawaii, 96822\label{fnm:1}.}%
% \footnotetext{MAPMO-F\'ed\'eration Denis Poisson, University of Orl\'eans, 45067 Orl\'eans, France\label{fnm:2}.}%

\begin{abstract}
Missions to rendezvous with or capture an asteroid present significant interest both from a geophysical and safety point of view. They are key to the understanding of our solar system are as well stepping stones for interplanetary human flight. In this paper, we focus on a rendezvous mission with \RH, an asteroid classified as a Temporarily Captured Orbiter (TCO). TCOs form a new population of near Earth objects presenting many advantages toward that goal. Prior to the mission, we consider the spacecraft hibernating on a Halo orbit around the Earth-Moon's $L_2$ libration point. The objective is to design a transfer for the spacecraft from the parking orbit to rendezvous with \RH while minimizing the fuel consumption. Our transfers use indirect methods, based on the Pontryagin Maximum Principle, combined with continuation techniques and a direct method to address the sensitivity of the initialization. We demonstrate that a rendezvous mission with \RH can be accomplished with low delta-v. This exploratory work can be seen as a first step to identify good candidates for a rendezvous on a given TCO trajectory.
\end{abstract}

% \begin{keyword}


% \MSC[2010] 49M05 \sep 70F07 \sep 93C10

% \end{keyword}

% \end{frontmatter}


\section{Introduction}
\label{intro}
The motivation for our work is to study asteroid capture missions for a specific population of near Earth objects. The targets, Temporarily Captured Orbiters (TCO), are small asteroids that become temporarily-captured on geocentric orbits in the Earth-Moon system.  They are characterized as satisfying the following constraints:
\begin{itemize}
\item the planetocentric Keplerian energy $E_{planet}<0$;
\item the planetocentric distance is less than three Earth's Hill radii (e.g., $3R_{H,\oplus} \sim 0.03$ AU);
\item it makes at least one full revolution around the Earth in the Earth-Sun co-rotating frame, while satisfying the first two constraints.
\end{itemize}

In regard to the design of a round trip mission, the main advantage of the TCOs lies in the fact that those objects have been naturally redirected to orbit the Earth. This contrasts with recently proposed scenarios to design, for instance, a robotic capture mission for a small near-Earth asteroid and redirect it to a stable orbit in the EM-system, to allow for astronaut visits and exploration (e.g. the Asteroid Redirect Mission (ARM)). 

\begin{figure}[ht!]
\centering
\includegraphics[width=6cm]{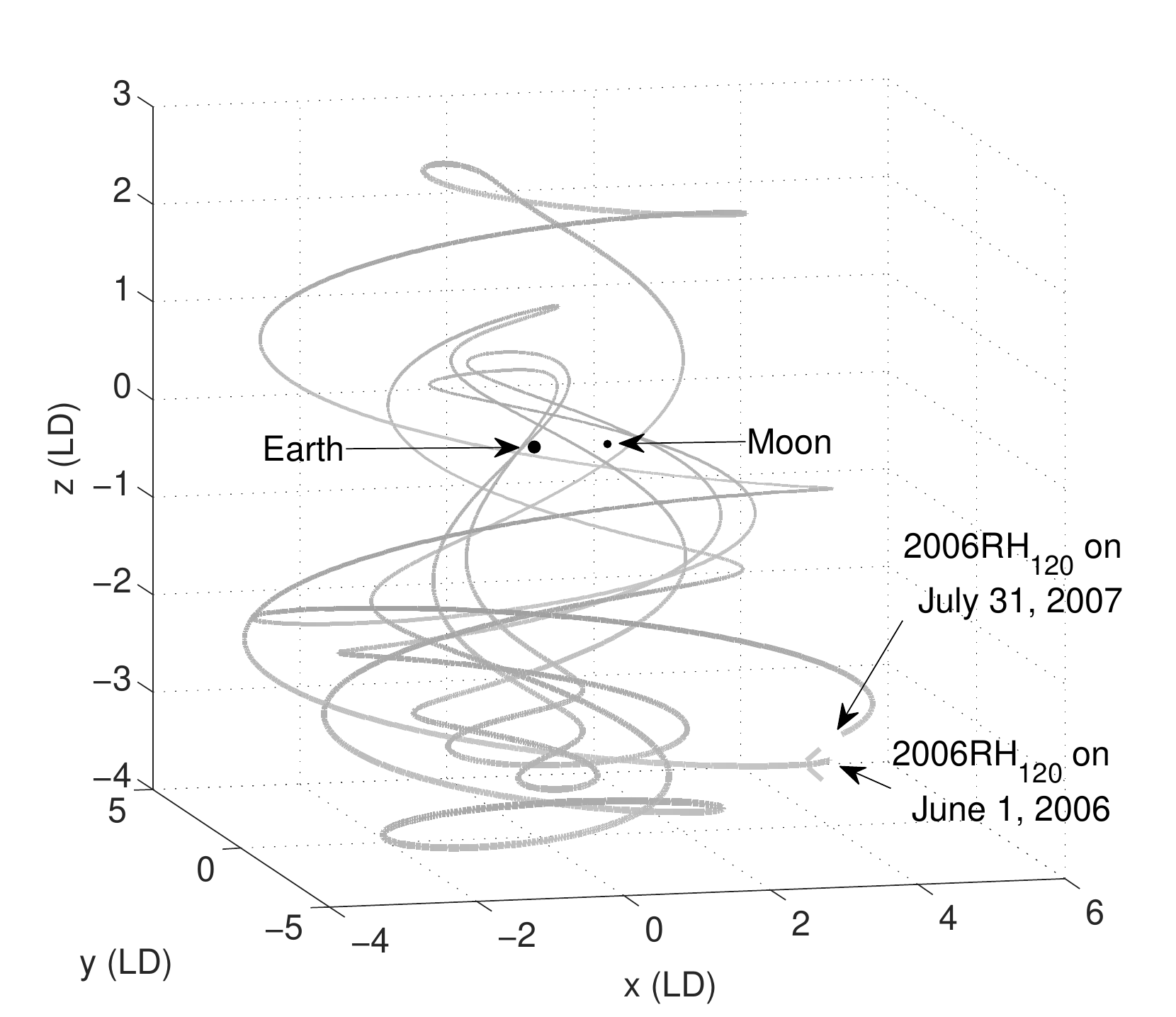}
\includegraphics[width=6cm]{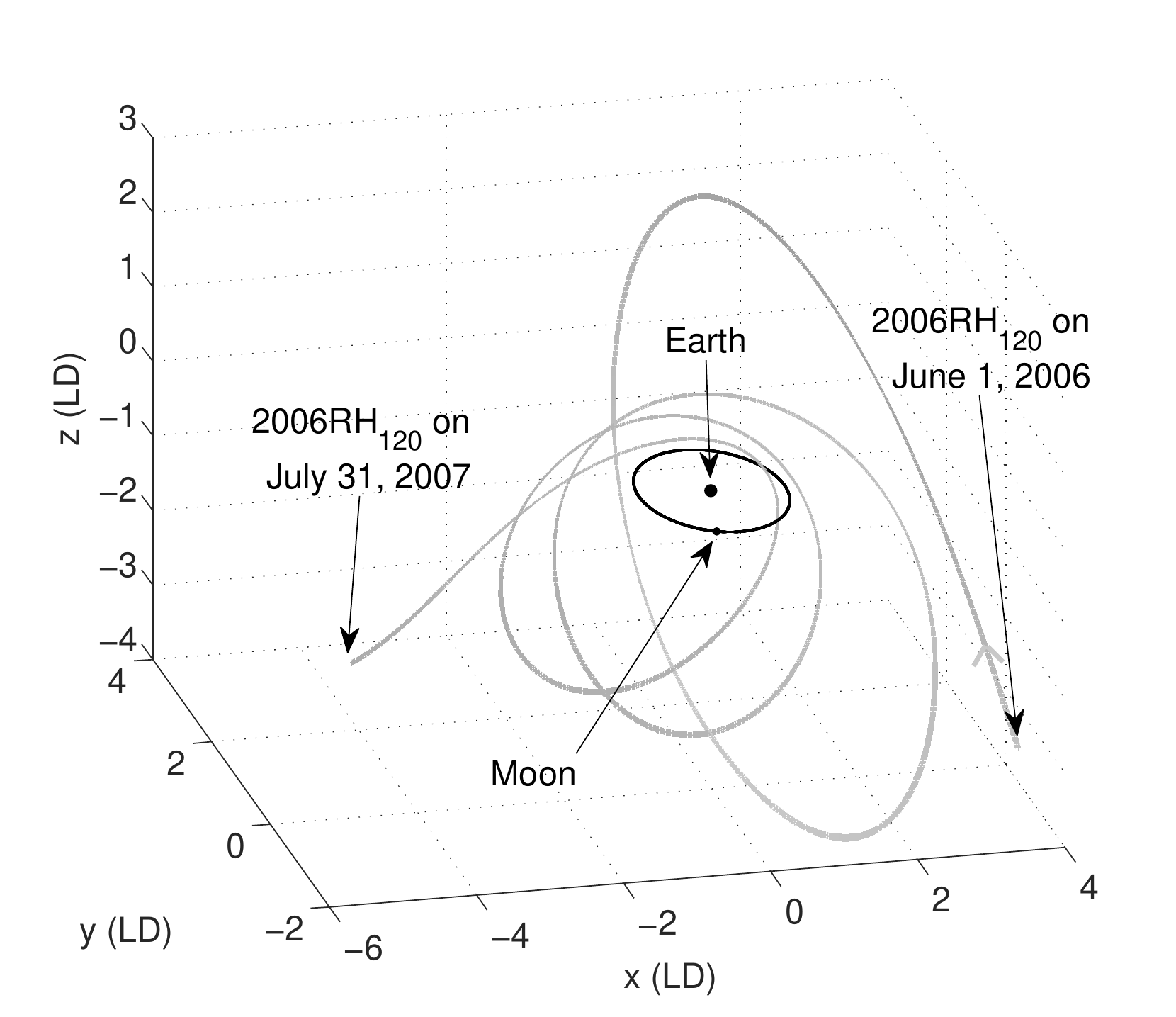}
\caption{Orbit of \RH in the Earth-Moon CR3BP rotating (left) and inertial (right) reference frame. \label{RH120fulltraj}}
%\qquad
\end{figure}

In this paper we focus on \RH, a few meters diameter near Earth asteroid, officialy classified as a TCO. \RH was discovered by the Catalina Sky Survey on September 2006. Its orbit from June 1st 2006 to July 31st 2007 can be seen on Figure \ref{RH120fulltraj}, generated using the Jet Propulsion Laboratory's HORIZONS database which gives ephemerides for solar-system bodies. The period June 2006 to July 2007 was chosen to include the portion of the orbit within the Earth's Hill sphere with a margin of about one month. We can also observe that \RH comes as close as 0.72 normalized units from Earth-Moon barycenter.
In \cite{Jedicke}, the authors investigate a population statistic for TCOs. Their work is centered on the integration of the trajectories for 10 million test-particles in space, in order to classify which of those become temporarily-captured by the Earth's gravitational field -- over eighteen-thousand of which do so. This suggests that \RH is not the only TCO and that it is relevant to compute a rendezvous mission to this specific asteroid to get insight whether TCOs can be regarded as possible targets for transfers with small fuel consumption, and thus cost.

In this paper, our goal is to solve the optimal control problem that consists in performing a rendezvous mission from a Halo orbit around the Earth-Moon libration point $L_2$ to the \RH asteroid, while minimizing the fuel consumption. The choice of targets for our rendezvous mission sets us apart from the existing literature where transfers are typically designed between elliptic orbits in the Earth-Moon or other systems \cite{Caillau,Mingotti1} and \cite{Mingotti2}, or to a Libration point \cite{Folta,Picot} and \cite{Ozimek,Vaquero}. Rendezvous missions to asteroids in the inner solar system can be found in \cite{Dunham,Hayabusa} but they concern asteroids on elliptic orbits which is not the case for us since TCOs are presenting complex orbits and therefore require a different methodology. Our assumption on the hibernating location for the spacecraft, a Halo orbit around the Earth-Moon unstable Libration points $L_2$, is motivated in part from the successful Artemis mission \cite{Arthemis1,Arthemis2} and in part from the constraint on the duration of the mission, mostly impacted by the time of detection of the asteroid. Indeed, the Arthemis mission demonstrated low delta-v ($\Delta v$) station keeping on Halo orbits around $L_1$ and $L_2$. While this first study focuses on similar Halo orbits to the Artemis mission, the study needs to be extended to analyze a variation of Halo orbits both in shape (to include the eight-shaped Lissajous orbits) and also in $z$-excursion or energy. As a first approach, in this work we search the best departure point on the prescribed Halo orbit without practical consideration related to synchronization with the asteroid orbit. The primary objective for the mission is to maximize the final mass, but since the Earth capture duration of the asteroid is limited it imposes a time constraint which need to be addressed eventually. In this paper however, we consider rendezvous missions possibly over long period of time to determine what is the correlation between the required $\Delta v$ and the duration of the mission. A forthcoming paper will bridge our results to the TCOs detection aspect. 

Many methodologies have been developed over the past decades to design optimal transfers in various scenarios. Due to the complexity of the TCO orbits and the nature of the mission, techniques based on analytical solutions such as in \cite{Kluever1} for circular Earth orbits are not suitable and we use a numerical approach. A survey on numerical methods can be found in \cite{Conway1}, and for reasons related to the specifics of our problem we choose to use a deterministic approach based on tools from geometric optimal control versus an heuristic method such as in \cite{Besette,Conway2} and \cite{Vaquero,Zhu}. More precisely, we use an indirect method based on the maximum principle, as well as a direct method and continuation techniques to address the difficulties of initialization for our numerical scheme. Additionally, we fix the structure of the  control norm to be a piecewise constant function with the magnitude of the thrust either zero or maximum with at most three maximum magnitude arcs. Our techniques are illustrated on the \RH TCO and four others coming from the database of \cite{Jedicke}. Validation of our approach can be seen by comparing the work to \cite{Dunham}, in which the authors develop a low $\Delta v$ asteroid rendezvous mission that makes use of a Halo orbit around Earth-Moon $L_2$. Their situation is, however, different from ours, in that they have carefully chosen a specific asteroid for rendezvous. With a one-year transfer time, the $\Delta v$ value they realize is 432 m/s, which is comparable to the $\Delta v$ values presented here in a less-ideal scenario.

\section{Optimal Control Problem and Numerical Algorithm}
\label{opti}
\subsection{Model}
\label{section-circular}
In this paper, the circular restricted three-body problem \cite{KLMR} is used to approximate the spacecraft dynamics. This is justified by the fact that a TCO can be assumed of negligible mass, and that the spacecraft evolving in the TCO's temporary capture space is therefore attracted mainly by two primary bodies, the Earth and the Moon.

The CR3BP model is well known, let us however go through some notations fo the sequel of the paper. We denote by $(x(t),y(t),z(t))$ the spatial position of the spacecraft at time $t$. In the rotating coordinates system, and under proper normalization, the primary planet identified here to the Earth, has mass $m_1=1-\mu$ and is located at the point $(-\mu,0,0)$; while the second primary, identified to the Moon, has a mass of $m_2=\mu$ and is located at $(1-\mu,0,0)$, where $\mu=1.2153e-2$. The distances of the spacecraft with respect to the two primary bodies are given by $\varrho_{1}=\sqrt{(x+\mu)^{2}+y^{2}+z^{2}}$, $\varrho_{2}=\sqrt{(x-1+\mu)^{2}+y^{2}+z^{2}}$ respectively. The potential and kinetic energy, respectively $V$ and $K$ of the system are given by 
\vspace{-0.1in}\begin{equation}
\label{eq:potential}
V=\frac{x^2+y^2}{2}+\frac{1-\mu}{\varrho_{1}}+\frac{\mu}{\varrho_{2}}+\frac{\mu(1-\mu)}{2},\;\;K=\frac{1}{2}(\dot x^2+\dot y^2+\dot z^2).
\end{equation}
We assume a propulsion system for the spacecraft is modeled by adding terms to the equations of motion depending on the thrust magnitude and some parameters related to the spacecraft design. The mass of the spacecraft is denoted by $m$.  Under those assumptions, we have the following equations of motion: 
\begin{eqnarray}\label{eq:motion}
\ddot{x} - 2\dot{y} = \frac{\partial V}{\partial x}+\frac{T_{\max}}{m}u_1, \quad
\ddot{y} + 2\dot{x} = \frac{\partial V}{\partial y}+\frac{T_{\max}}{m}u_2, \quad
\ddot{z}= \frac{\partial V}{\partial z}+\frac{T_{\max}}{m}u_3 
\end{eqnarray}
where $u(.)=(u_1(.),u_2(.),u_3(.))$ is the control, and satisfies the constraint $\|u\|=\sqrt{u_1^2+u_2^2+u_3^2}\leq 1$, and with $T_{\rm max}/m$ normalized. A first integral of the free motion is given by the energy of the system $E=K-V$. We will later use this energy value to analyze the choice of the rendezvous point and the parking orbit for the spacecraft. It is well know that the uncontrolled motion of the dynamical system has five equilibrium points defined as the critical points of the potential $V$. Three of them $L_1,L_2$ and $L_3$ are aligned with the Earth-Moon axis and have been shown to be unstable, while the two others are stable and are positioned to form equilateral triangles in the plane of orbit with the two primaries. Since our goal is to maximize the final mass we must include the differential equation governing the variation of the mass along the transfer:
\begin{equation}
\dot m = -\beta T_{\max} \|u\|
\end{equation}
 where the parameter $\beta$, the thruster characteristic of our spacecraft, is given by $\beta=\frac{1}{I_{sp}g_0}$ (it is the inverse of the ejection velocity $v_e$), with $I_\mathrm{sp}$ the specific impulse of the thruster and $g_0$ the acceleration of gravity at Earth sea level.
 %%%
 %%%
 %%%
\subsection{Optimal Control Problem}
Let $q=(q_s,q_v)^t$ where $q_s=(x,y,z)^t$ represents the position variables and $q_v=(\dot x,\dot y,\dot z)^t$ the velocity ones. From section \ref{section-circular}, our dynamical system is an affine control system of the form:
\begin{equation}
\label{eq-affine}
 \dot{q}=F_0(q)+\frac{T_{\max}}{m}\sum_{i=1}^3F_i(q)u_i
\end{equation}
where the drift is given by $F_0(q)=(\dot x,\dot y,\dot z,2\dot y+x-\frac{(1-\mu)(x+\mu)}{\varrho_1^3}-\frac{\mu(x-1+\mu)}{\varrho_2^3},
-2\dot x+y-\frac{(1-\mu)y}{\varrho_1^3}-\frac{\mu y}{\varrho_2^3},-\frac{(1-\mu)z}{\varrho_1^3}-\frac{\mu z}{\varrho_2^3})^t$, and 
 the control vector fields are $F_i(q)=\vec e_{3+i}$ with $\vec e_i$ forming the orthonormal basis of $\R^6$. 

We consider the rendezvous transfer from an initial point $q^{rdv}(t_0)$ on a parking orbit $\mathcal{O}_0 \in \R^6$ to a final position and velocity $q^{rdv}(t_f)$ on the \RH orbit. Note that the initial and final positions and velocities are variables of the global optimization problem. The criterion to maximize is the final mass which is equivalent to minimizing the fuel consumption or the $\Delta v=\int_{t_0}^{t_f} \frac{T_{\max}\|u(t)\|}{m(t)}dt$. Since the mass evolves proportionally to the norm of the thrust, our criterion is equivalent to the minimization of the $L_1$-norm of the control:
\begin{eqnarray}\label{opt_cont_pb_3d}
   \min_{u\in{\mathcal U}} \int_{t_{0}}^{t_{f}} \|u(t)\|dt,
\end{eqnarray}
where $\mathcal U=\{u(.);{\rm measurable\;bounded\;and}\;\|u(t)\|\leq 1\;{\rm for\;almost \;all}\;t\} $, $t_0$ and $t_f$ are respectively the initial and final time. Remark that since we choose $\mathcal{O}_0$ to be a Halo orbit around a libration point, it is uniquely determined by a single point of the orbit using the uncontrolled CR3BP dynamics. This fact will play an important role for one of the necessary optimality conditions below. The large number of variables in our problem significantly adds complexity to the search for a solution. In particular, in the case of free final time we expect an infinite time horizon with a control structure that mimics impulsive maneuvers. To simplify our optimal problem we have two options, either we fix the transfer time or we fix the structure of the control. If we fix the final time, the sensitivity of the shooting method can be addressed by using the solution of a smoother criterion than the $L_1$-norm, for instance the $L_2$-norm, and linking it to the target criterion by a continuation procedure, see \cite{Haberkorn}. In this paper, however, we take a different approach and decide to fix the structure of the control. In the sequel we focus on designing transfers associated to controls with a piecewise constant norm with value in $\{0,1\}$ and four switchings. More precisely, we consider control functions that are piecewise continuous such that there exists times 
$t_0 \leq t_1 \leq t_2 \leq t_3 \leq t_4 \leq t_f$ with
\begin{equation}
\|u(t)\| =\left\{\begin{array}{cc}
1 & \ {\rm if}\ t \in (t_0;t_1) \cup (t_2;t_3) \cup (t_4;t_f)\\
0 &  \ {\rm if}\ \in (t_1;t_2) \cup (t_3;t_4).
\end{array}\right.
\end{equation}
\begin{figure}[ht!]
\centering
\includegraphics[width=7cm]{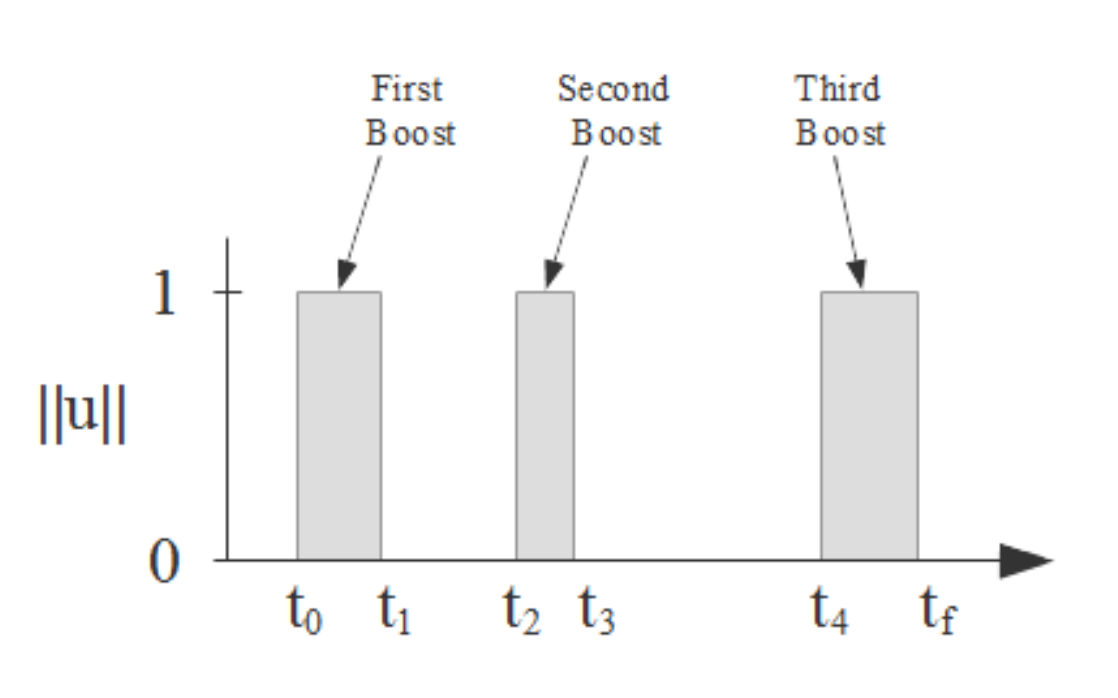}
\caption{Control function modeling thrust impulses over time. \label{thrustFunction}}
\end{figure}
Here the final time $t_f$ is free. Note that  our numerical method will be able to select the best control strategy even if it has less than three boosts.
\subsection{Necessary Conditions for Optimality}\label{s_CN}

The maximum principle provides first order necessary conditions for a trajectory to be optimal \cite{Pont}. Details regarding the application of the maximum principle to orbital transfers can be found in many references, including \cite{Caillau,Pont}. We denote by $X(t) = (q(t),m(t)) \in \R^{6+1}$ the state, where  $q=(x,y,z,\dot{x},\dot{y},\dot{z})$ is the position and velocity of the spacecraft with $m$ its mass. The maximum principle applied to our optimal control problem states that if $(q(\cdot),m(\cdot),u(\cdot)) = (X(\cdot),u(\cdot))$ is an optimal solution defined on $[t_0,t_f]$, then there exists an absolutely continuous adjoint state $(p^0,p_X(\cdot)) = (p^0,p_q(\cdot),p_m(\cdot))$, defined on $[t_0,t_f]$ such that:
\begin{itemize}
\item
	$(p^0,p_X(\cdot)) \neq 0,\ \forall t \in [t_0,t_f]$, and $p^0 \leq 0$ is a constant.
\item
	Let $H$, the Hamiltonian, be $H(t,X(t),p^0,p_X(t),u(t)) = p^0 \|u(t)\| + \langle p_X(t) , \dot{X}(t)\rangle$, then
	\begin{eqnarray}\label{Hdyn1}
	\dot{X}(t) = \frac{\partial H}{\partial p_X}(t,X(t),p^0,p_X(t),u(t)),\ {\rm for\ a.e.}\ t \in [t_0,t_f],\\ \label{Hdyn2}
	\dot{p}_X(t) = -\frac{\partial H}{\partial X}(t,X(t),p^0,p_X(t),u(t)),\ {\rm for\ a.e.}\ t \in [t_0,t_f],
	\end{eqnarray}
	where $\langle ,\rangle$ denotes the inner product. 
\item
	$H(t,X(t),p^0,p_X(t),u(t)) = \max_{\|\nu\| \leq 1} H(t,X(t),p^0,p_X(t),\nu),\ \forall t$ s.t. $\|u(t)\| = 1$ (maximization condition).
\item
	$\Psi(t_{i}) = 0$ for $i=1,\cdots,4$.
\item
	$H(t_f,X(t_f),p^0,p_X(t_f),u(t_f)) = 0$, if $t_f$ is free.
\item
	$\langle p_q(t_0) ,F_0(q(t_0)) \rangle = 0$ (initial transversality condition).
\item
	$p_m(t_f) = 0$.
\end{itemize}
The function $\Psi(\cdot)$ is the so-called switching function corresponding to the problem with an unrestricted control strategy and we have $\Psi(t)=p_0 + T_\mathrm{max}\Big(\frac{\|p_v(t)\|}{m(t)}-p_m(t)\beta \Big)$.

The maximization condition of the Hamiltonian $H$ is used to compute the control on $[t_0;t_1]\cup[t_2;t_3]\cup[t_4;t_f]$ and we have $u(t) = \frac{p_v(t)}{\|p_v(t)\|}\ {\rm for\;all}\; t \in [t_0;t_1]\cup[t_2;t_3]\cup[t_4;t_f]$. The initial transversality condition reflects the fact that the initial departing point is free on the Halo orbit $\mathcal{O}_0$. Remark that since the data for the TCO's trajectory are given as ephemerides, there are no dynamics equations to describe those orbits in the CR3BP. Thus, we cannot compute the tangent space to a TCO point and we cannot extract a transversality condition for $p_q(t_f)$ at the rendezvous point. Since we expect numerous local extrema for this optimal control problem, it is however preferable to solve the problem for fixed rendezvous points on a discretization of the TCO orbit.

\subsection{Shooting Method}\label{s_shoot}
Our numerical method is based on the necessary conditions given in section \ref{s_CN}. Let $Z(t) = (X(t),p_X(t)),\ t \in [t_0;t_f]$, and $u(q,p)$ the feedback control expressed using the maximization condition. Then, we have $\dot{Z}(t) = \Phi(Z(t))$ where $\Phi$ comes from equations (\ref{Hdyn1}) and (\ref{Hdyn2}). The goal is to find $Z(t_0)$, $t_1$, $t_2$, $t_3$, $t_4$ and $t_f$ such that the following conditions are fulfilled:
\begin{enumerate}
\item
	$\Psi(t_i)=0$ for  $i=1,\cdots 4$;
\item
	$X(t_f)$ is the prescribed rendezvous point;
\item
	$X(t_0)\in \mathcal O_0$ and the initial transversality condition is verified;
\item
	$H(t_f)=0.$
\end{enumerate}
The problem has been transformed into solving a multiple points boundary value problem. More specifically, we must find the solution of a nonlinear equation $S(Z(t_0),t_1,t_2,t_3,t_4,t_f) = 0$, where $S$ is usually called the shooting function. 
The evaluation of the shooting function is performed using the high order numerical integrator \emph{DOP853}, see \cite{Hairer}. The search for a zero of the shooting function is done with the quasi-Newton solver \emph{HYBRD} of the Fortran \emph{minpack} package. Since $S(.)$ is nonlinear, the Newton-like method is very sensitive to the initial guess. This leads us to consider heuristic initialization procedures. We combine two types of techniques, a direct method and a continuation method. The discretization of the TCO's orbit requires the study of thousands of transfers, we use a direct method for a dozen of rendezvous and expand to other points on the orbit using a continuation scheme. The motivation is that direct methods are very robust but time consuming  while continuation method succeeded for our problem in most cases very rapidly. The direct method uses the modeling language \emph{Ampl}, see \cite{Ampl}, and the optimization solver \emph{IpOpt}, see \cite{IpOpt}, with a second order explicit Runge-Kutta scheme. Details on advanced continuation methods can be found in \cite{Haberkorn} and \cite{Hampath}, but we use a discrete continuation which is enough for our needs.

Moreover, note that in order to fulfill the initial transversality condition, we first prescribe $X(t_0)$ on the parking orbit and find a zero of the shooting function satisfying all the other conditions. Afterward, we do a continuation on $X(t_0)$ along the departing parking orbit in the direction that increases the final mass. Once a framing of the best $X(t_0)$ with respect to the final mass is found, we perform a final single shooting to satisfy the initial transversality condition along with the other conditions. This decoupling in the search of an extremal is motivated by the fact that we could very likely find a local maximum on $X(t_0)$ rather than a local minimum because of the periodicity of the initial parking orbit. We avoid this fact by first manually ensuring that the $X(t_0)$ we find will be the one for the best final mass and not the worst. However, our continuation procedure on $X(t_0)$ does not always succeed, mainly because of the high nonlinearity of the shooting function, as the trajectories we find can be very long. Even if some of the extremals we find do not satisfy the initial transversality condition with the aimed accuracy (typically a zero of the shooting function is deemed acceptable if $\|S(Z(t_0),t_{1,2,3,4,f})\| \leq 10^{-8}$), they are still rather close to satisfy it (of the order of $10^{-4}$).
%%%
%%%
%%%
\section{Results}
\subsection{Rendezvous to \RH}

\begin{figure}[ht!]
\includegraphics[width=14cm]{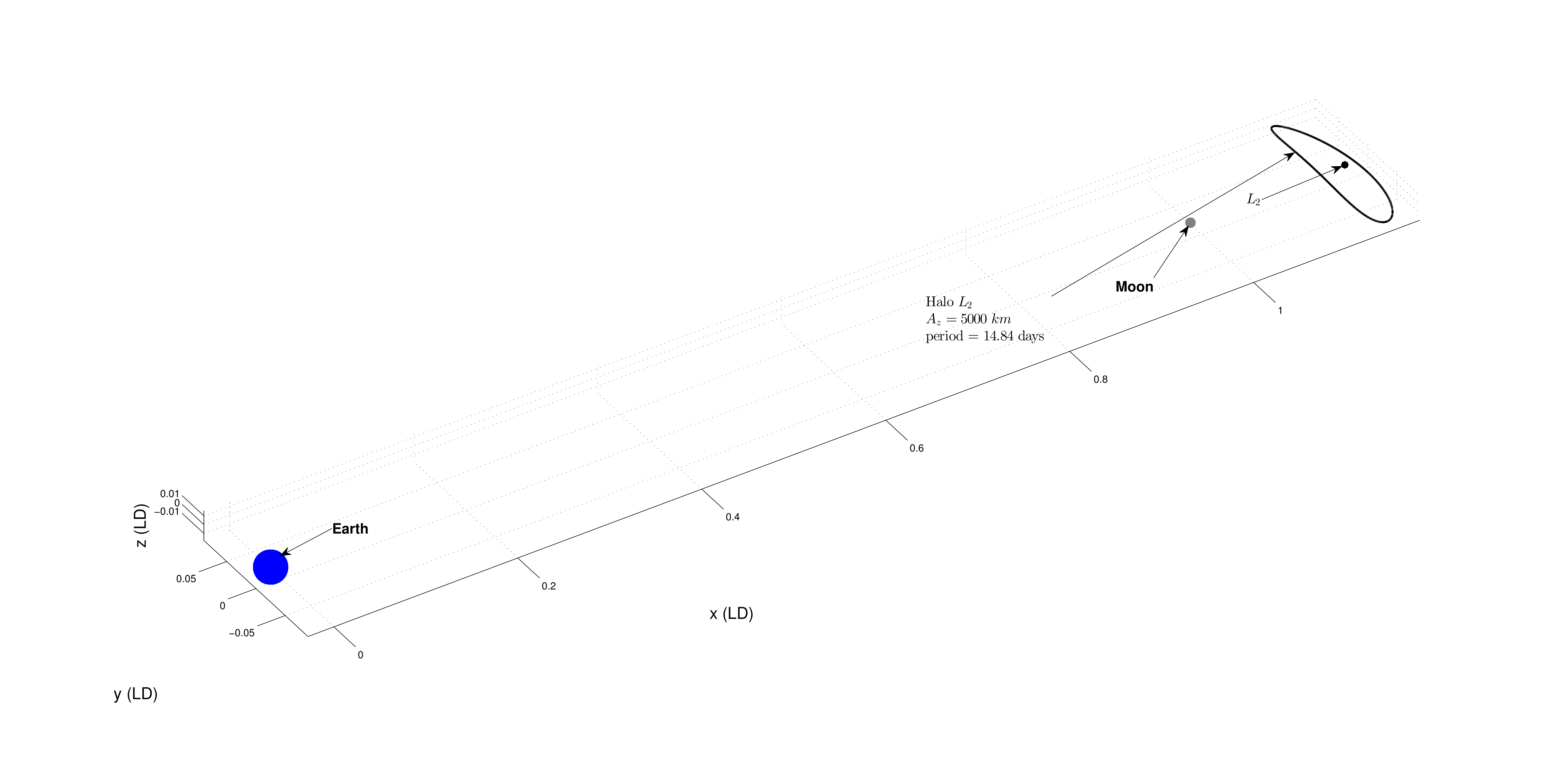}
\caption{Halo orbit from which the spacecraft is departing, $z-$excursion of 5000 km around the EM libration point $L_2$.}
\label{Halo}
\end{figure}

The objective of this section is to present an analysis of the evolution of the fuel consumption with respect to the location of the rendezvous point on the \RH orbit. As mentioned before, we restrict the study to transfers with at most three boosts. The primary goal is to obtain insights on the variation of fuel consumption based on the features of the rendezvous point such as its distance from $L_2$, or its energy. The spacrecraft characteristics are assumed to be an initial mass of $350\ kg$, a specific impulse $I_\mathrm{sp}$ of $230\ s.$ and a maximum thrust $T_\mathrm{max}$ of $22\ N$. The Halo orbit from which the spacecraft is departing is chosen to have a $z-$excursion of 5000 km around the EM libration point $L_2$, see Figure \ref{Halo}. The choice for this specific Halo orbit is motivated by the successful  Artemis mission \cite{Arthemis1,Arthemis2} that used Halo orbits with similar characteristics. The point corresponding to the positive $z-$excursion is $q_\mathrm{HaloL_2} = (1.119, 0, 0.013,$ $ 0., 0.180, 0)$, and the period of this particular Halo orbit is $t_\mathrm{HaloL_2} = 3.413$ in normalized time units or approximatively $14.84\ days$. 

\begin{figure}[ht!]
\includegraphics[width=6cm]{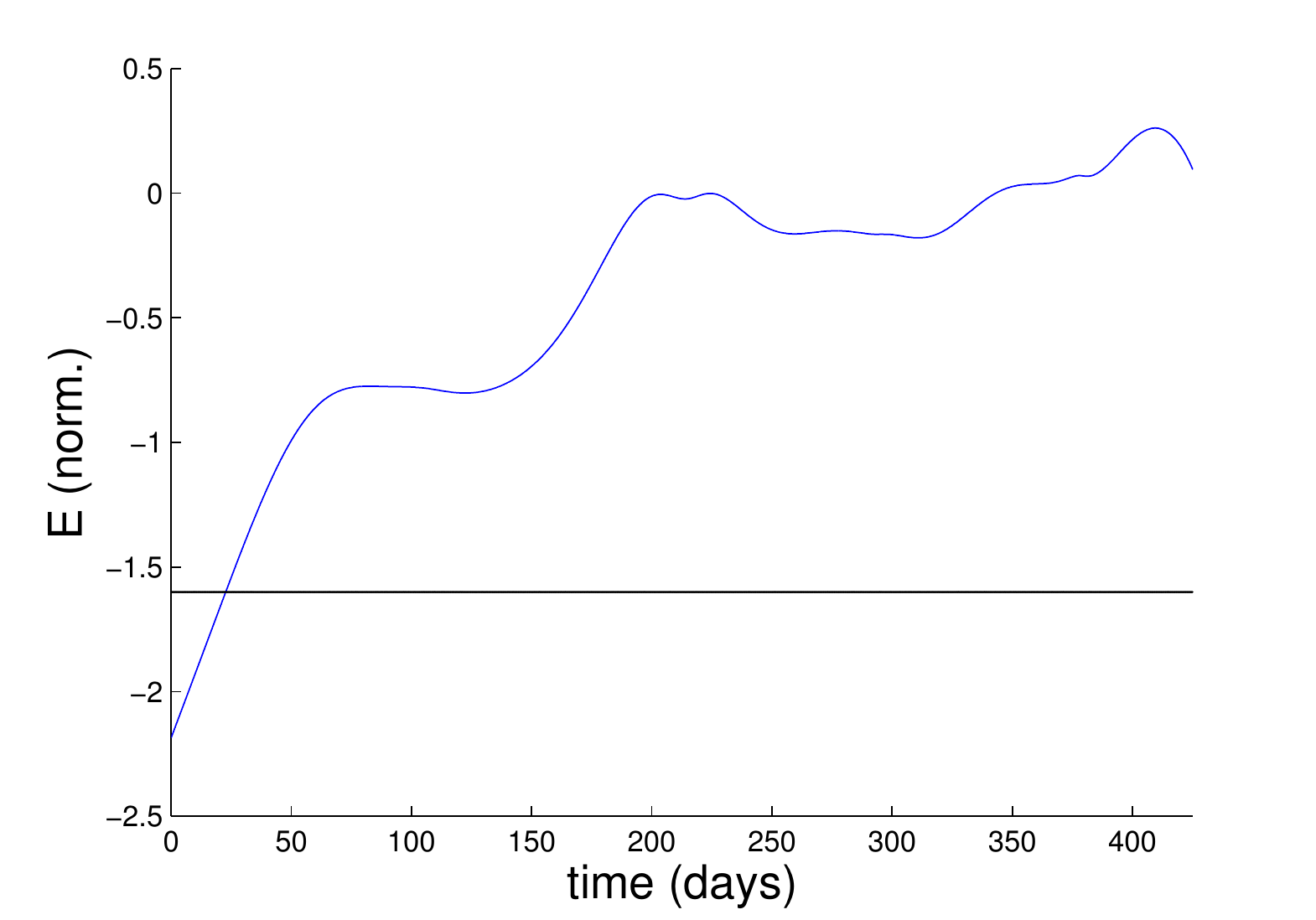} \includegraphics[width=8.3cm]{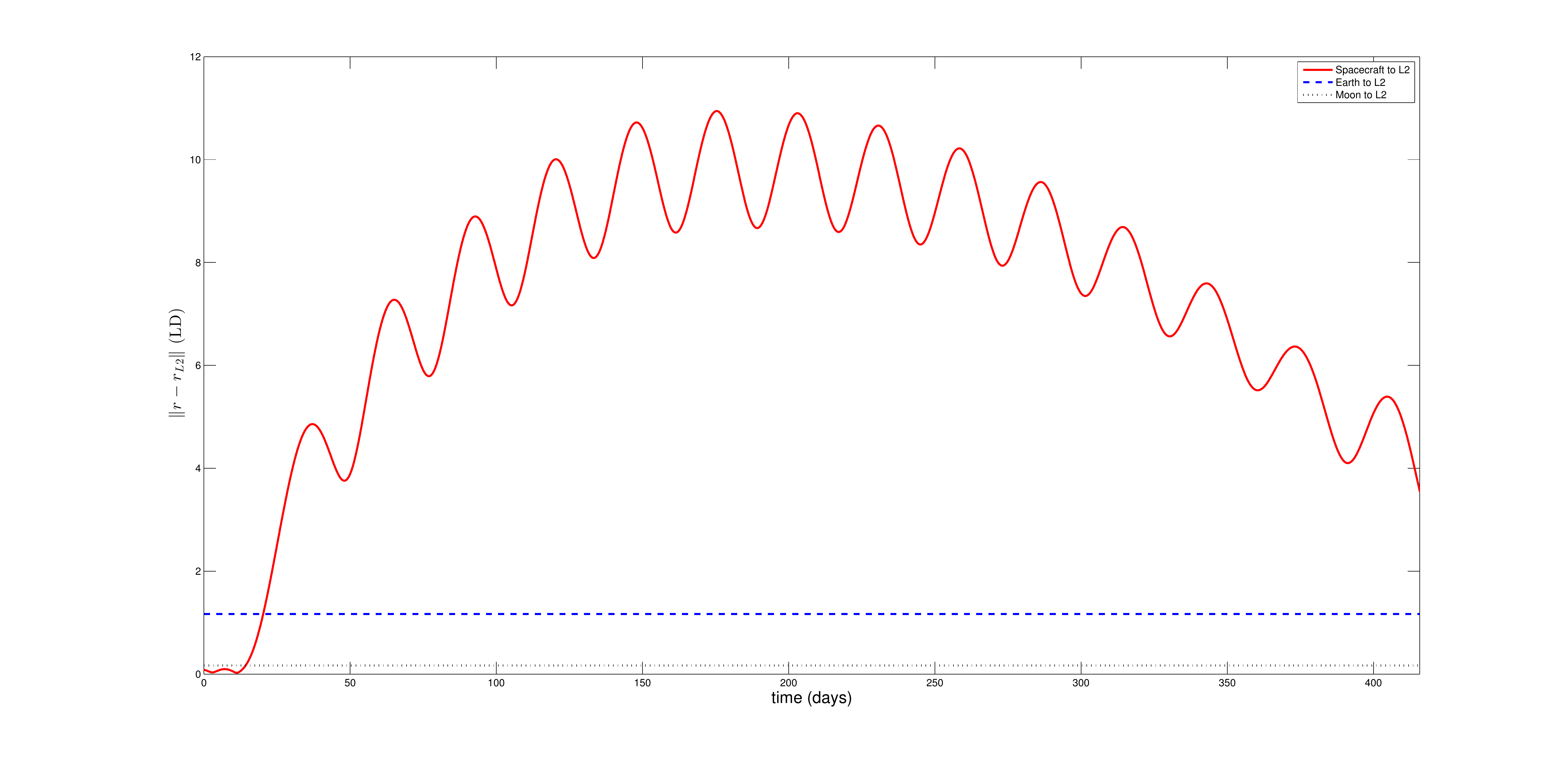}
\caption{Evolution of the energy (left) and of the Euclidean distance to the $L_2$ Libration point (right) for \RH. For energy, the  horizontal line represents the energy of the Halo periodic orbit around $L_2$ which is about -1.58.}
\label{RH120_Nrg_Dist}
\end{figure}

During the period represented in Figure \ref{RH120fulltraj}, asteroid  \RH does 17 clockwise revolutions around the origin of the CR3BP frame, and 3.6 revolutions in Earth inertial reference frame. The evolution of the energy of \RH and its distance to the Earth-Moon libration point $L_2$ as the asteroid evolves on its orbit are given in Figure \ref{RH120_Nrg_Dist}.

\begin{figure}[ht!]
\begin{center}
\includegraphics[width=6cm]{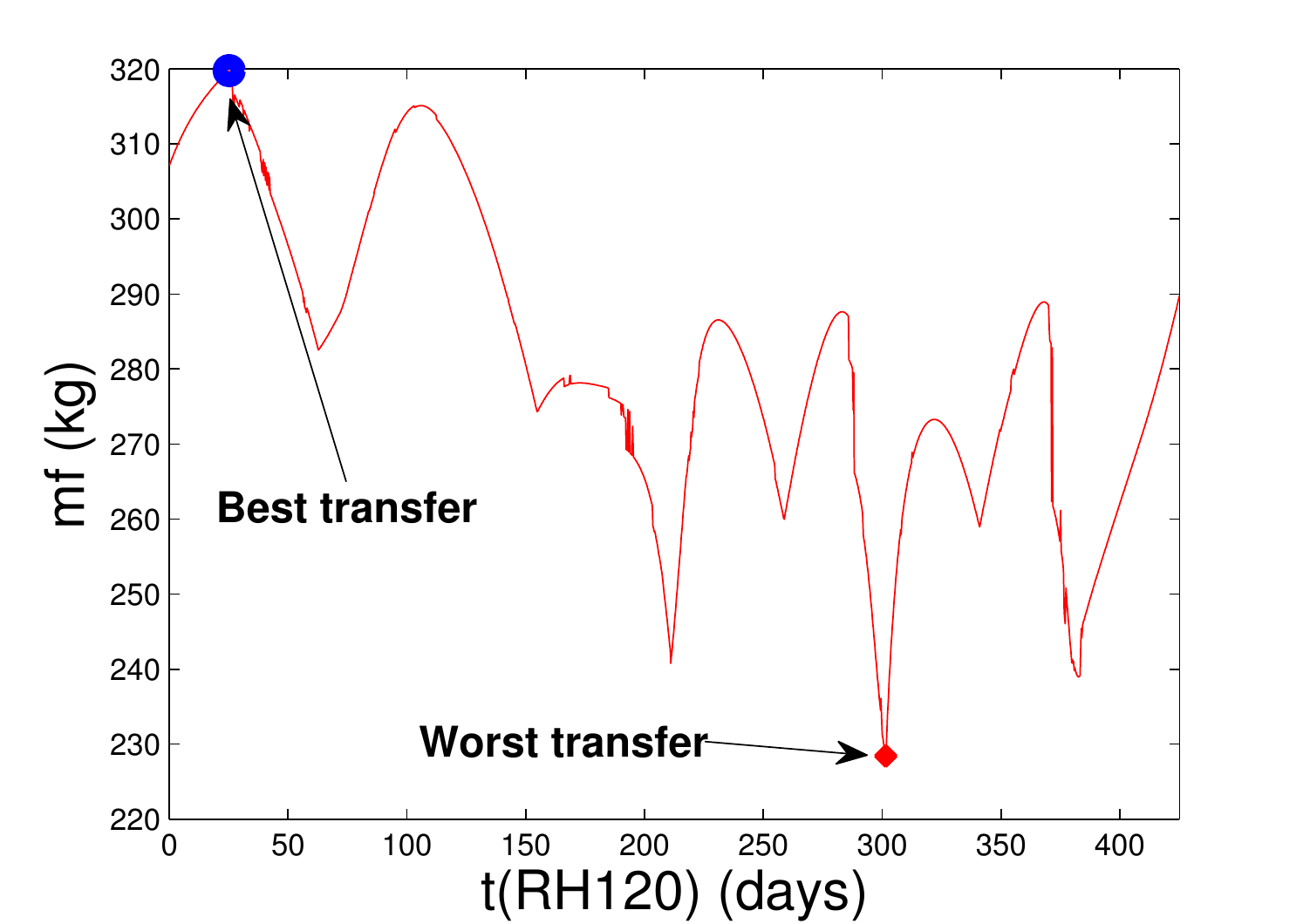} 
\includegraphics[width=6cm]{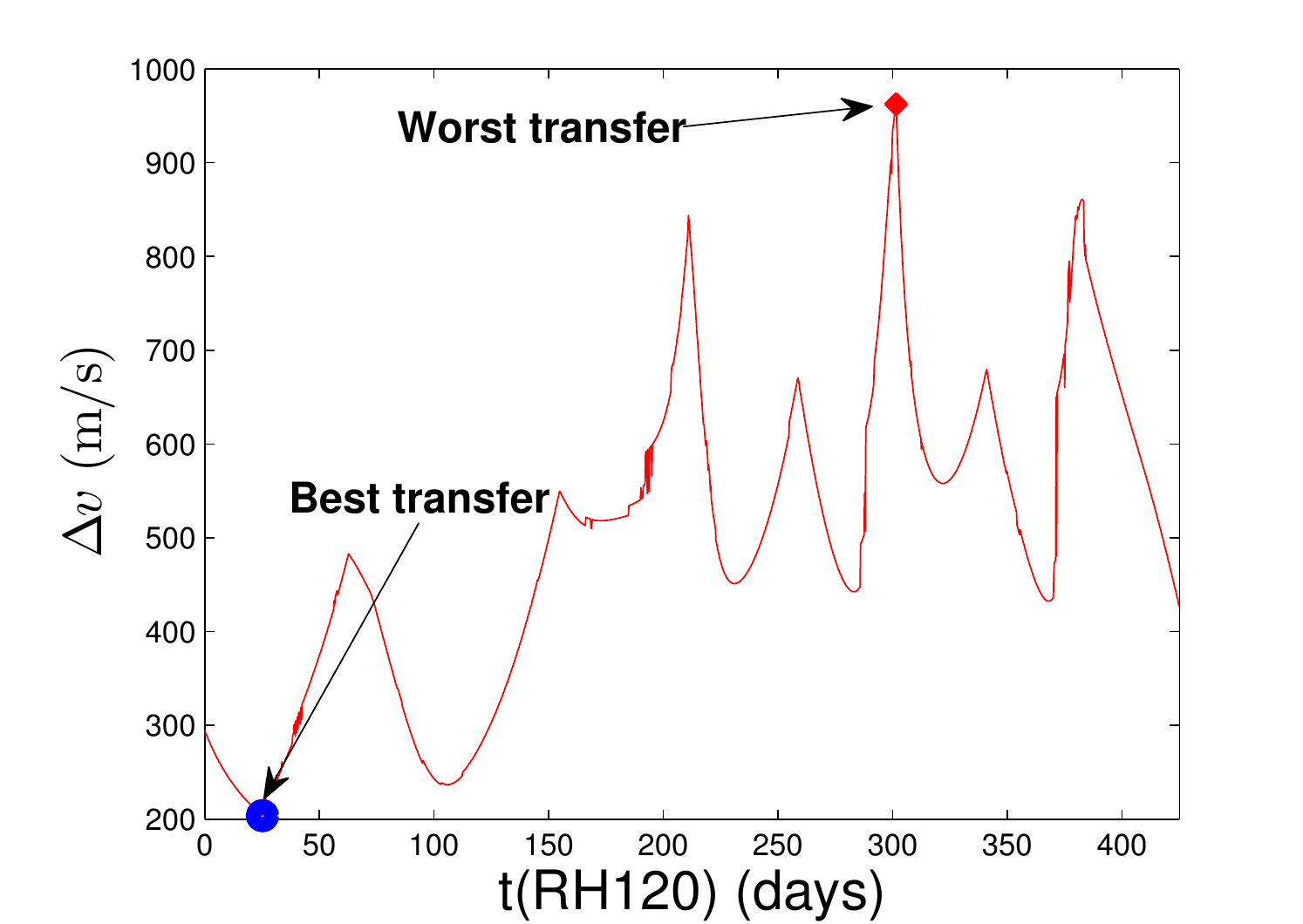}
\includegraphics[width=6cm]{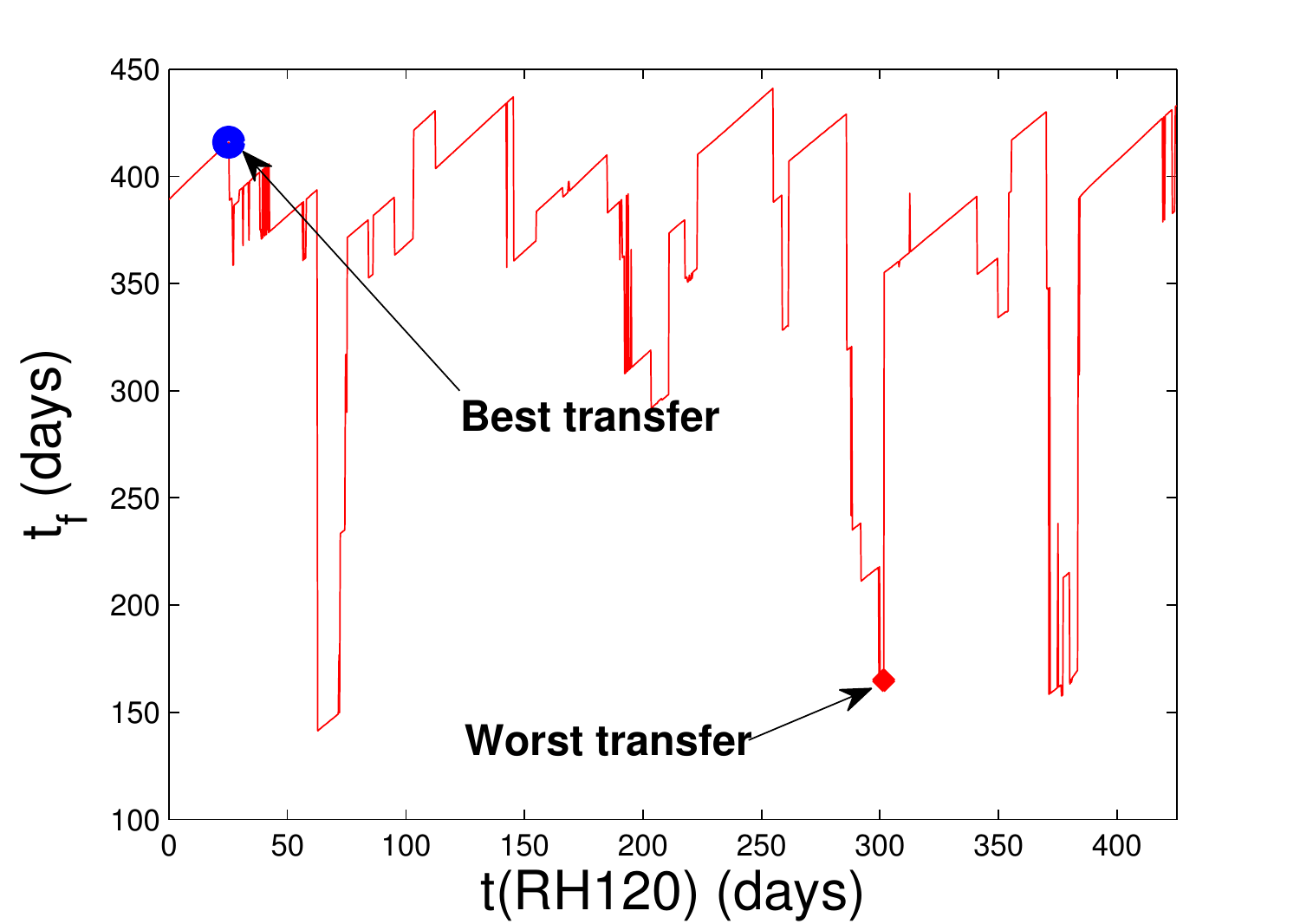}
\end{center}
\caption{\label{HaloL2_mfevo} Evolution of the final mass (top left), of the $\Delta v$ (top right) and transfer duration (bottom) from the Halo orbit around $L_2$ with respect to rendezvous points on the orbit of \RH.}
\end{figure}

\begin{figure}[ht!]
\begin{center}
\includegraphics[width=5cm]{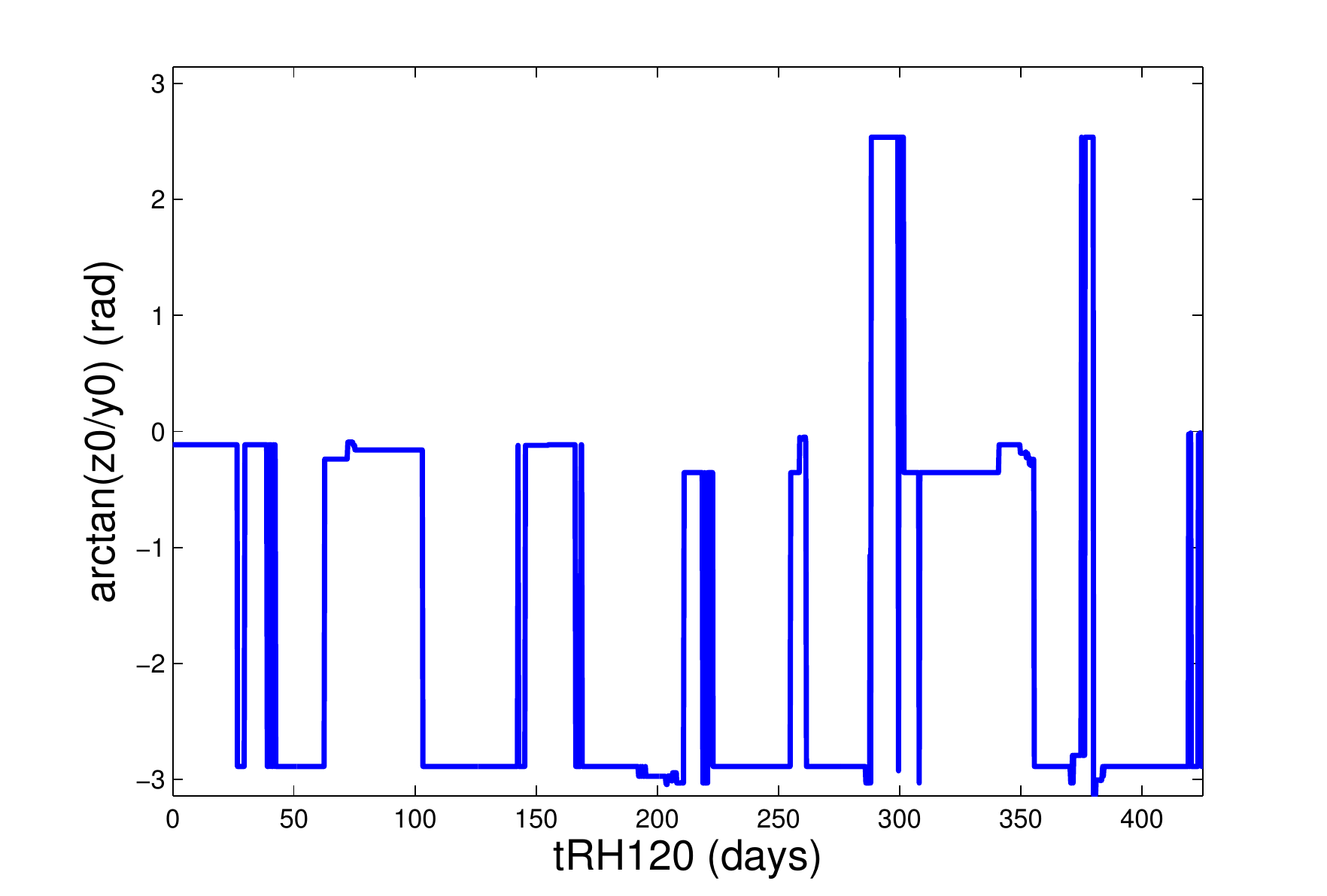}
\includegraphics[width=8cm]{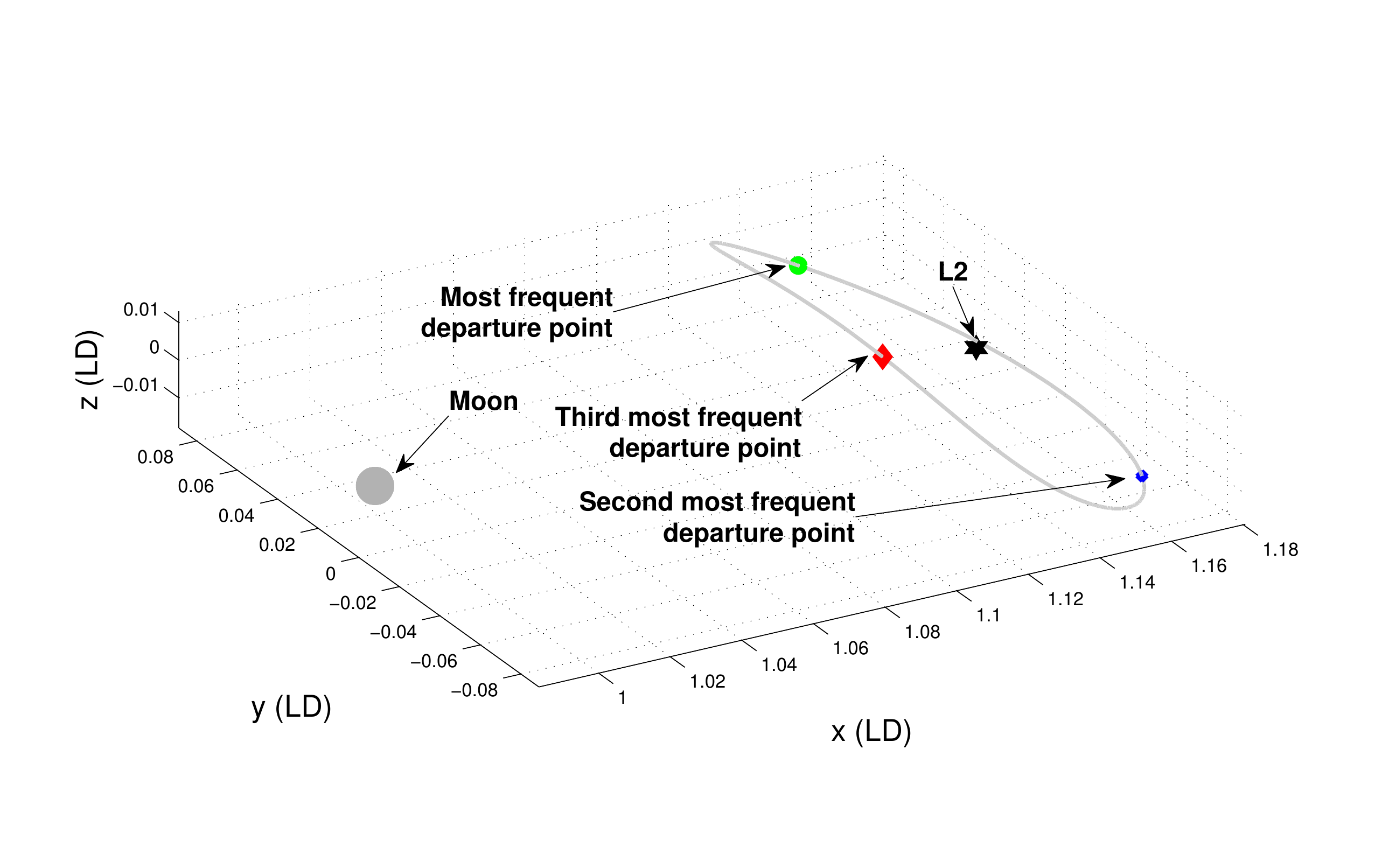}
\end{center}
\caption{\label{Angle_vs_rdv} The left picture represents the evolution of $(y,z)(0)$ argument with respect to time of rendezvous on \RH. The right picture represents the three most frequent departure points on the initial Halo orbit.}
\end{figure}

To analyze the variations of the fuel consumption with respect to the rendezvous point on the orbit, we discretize uniformly the orbit of \RH using 6 hours steps. For each rendezvous point of this discretization, we compute an extremal transfer (i.e. a solution of the maximum principle) with free final time using the techniques explained in section \ref{opti}. Figure \ref{HaloL2_mfevo} shows the evolution of the final mass, the $\Delta v$ and of the duration of the transfer with respect to the rendezvous point on the \RH orbit for a spacecraft departing from the Halo $L_2$ orbit and corresponding to the three boosts control strategy. As explained in section \ref{s_shoot}, some of the departure point on the Halo orbit are not fully optimized -- this is the case for about two thirds of them. Also note that a departure point different from $q_\mathrm{HaloL_2}$ implies a drift phase whose duration is not included in the transfer duration. It can be observed that the final mass has many local extrema and that the variation of the duration of the mission is not continuous (contrary to what we would expect). It is most likely due to local minima or to the fact that the value function of the problem can present discontinuities. 
Figure \ref{Angle_vs_rdv} shows the evolution of the departure point on the hibernating orbit for the spacecraft with respect to the rendezvous point on \RH (left) as well as the three most frequent departure points on the initial Halo orbit (right). More precisely, we represents the optimized argument of $(y,z)(0)$ (up to the quadrant: $arctan(z(0)/y(0))$). Note that the initial position on the Halo orbit directly depends on the optimized initial drift time. Since this initial drift phase has not always been successfully optimized, this figure has to be interpreted with caution. However, we can see that the departure position on the initial periodic orbit seems to always be close to a multiple of $\pi$.

Comparing the evolution of the final mass from Figure \ref{HaloL2_mfevo} and the evolution of \RH energy from Figure \ref{RH120_Nrg_Dist}, we can see that the best final masses are obtained on the first half of \RH trajectory, that is for rendezvous points which energies are closer to the departing energy. For the best transfer, the energy difference between the rendezvous point on \RH and the departing point on the Halo orbit is about 0.046. The rendezvous point with the closest energy to the initial orbit is only slightly before the optimal rendezvous point. Its final mass is 319.3 kg which is only 0.5 kg worst than the best final mass. This remark suggests strongly that a small difference in energy between the rendezvous point on \RH orbit and the departing point for the spacecraft on the Halo orbit is advantageous. 

\begin{figure}[ht!]
\begin{center}
\includegraphics[width=7cm]{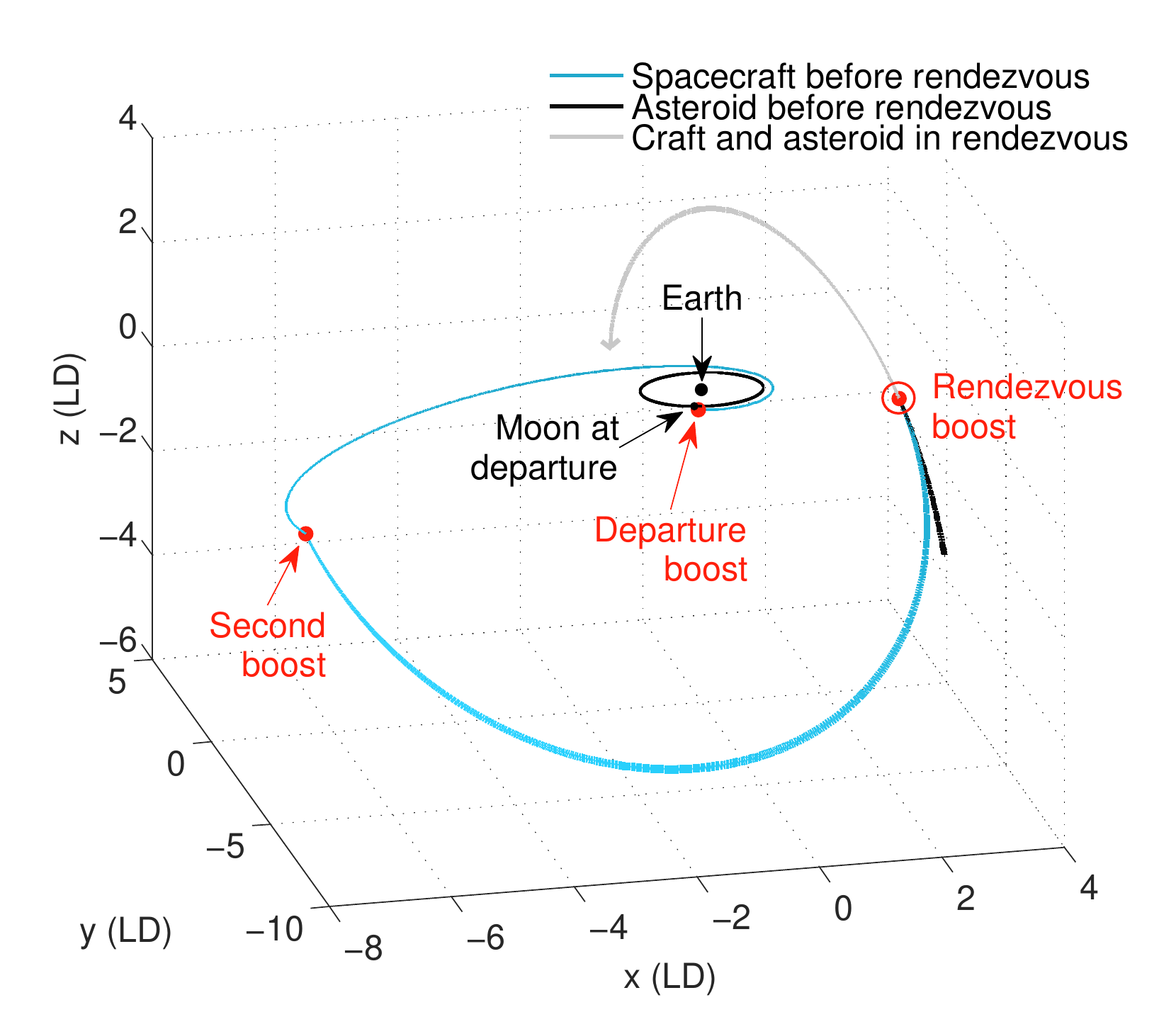}\includegraphics[width=7cm]{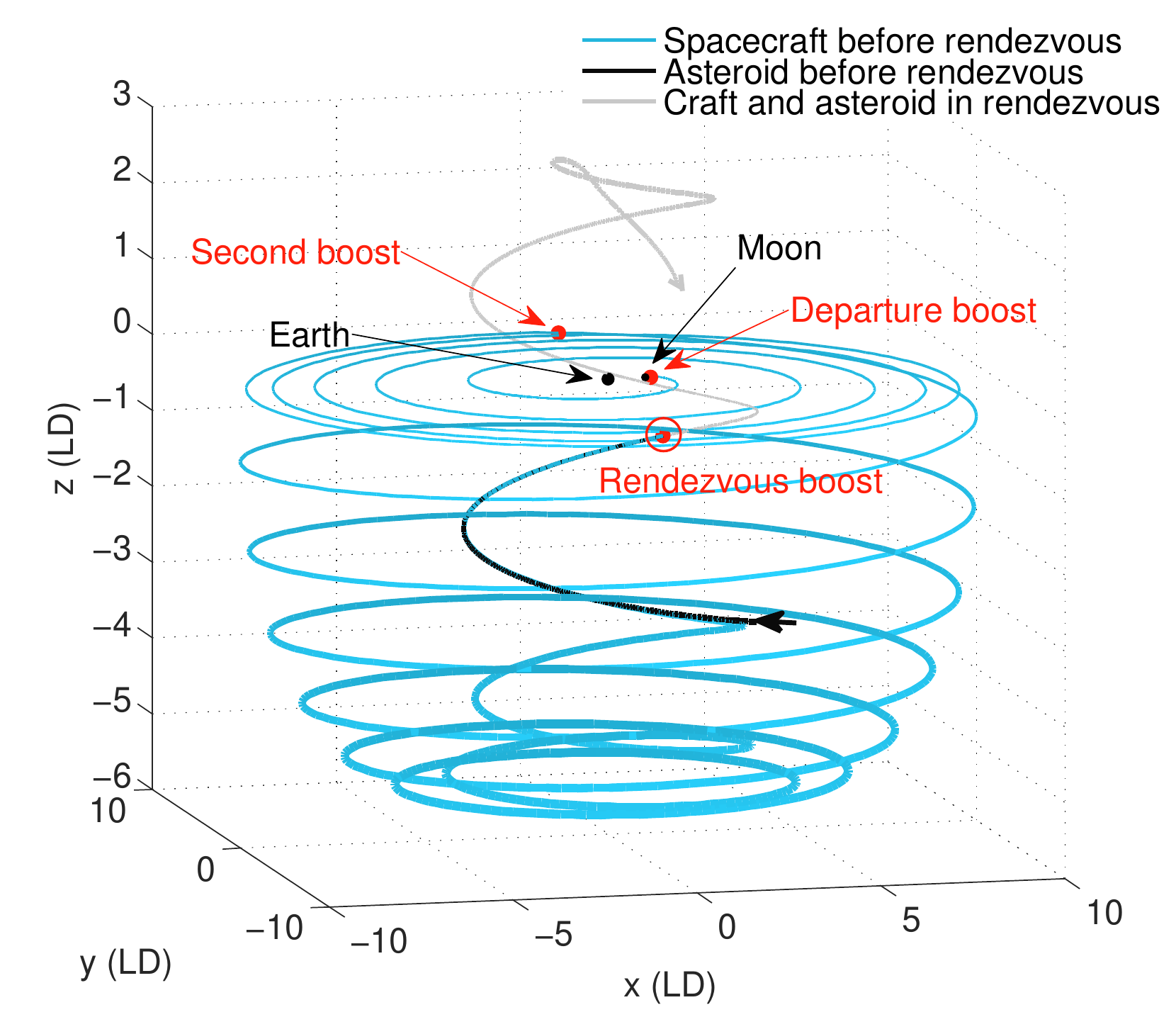}\\\includegraphics[width=9cm]{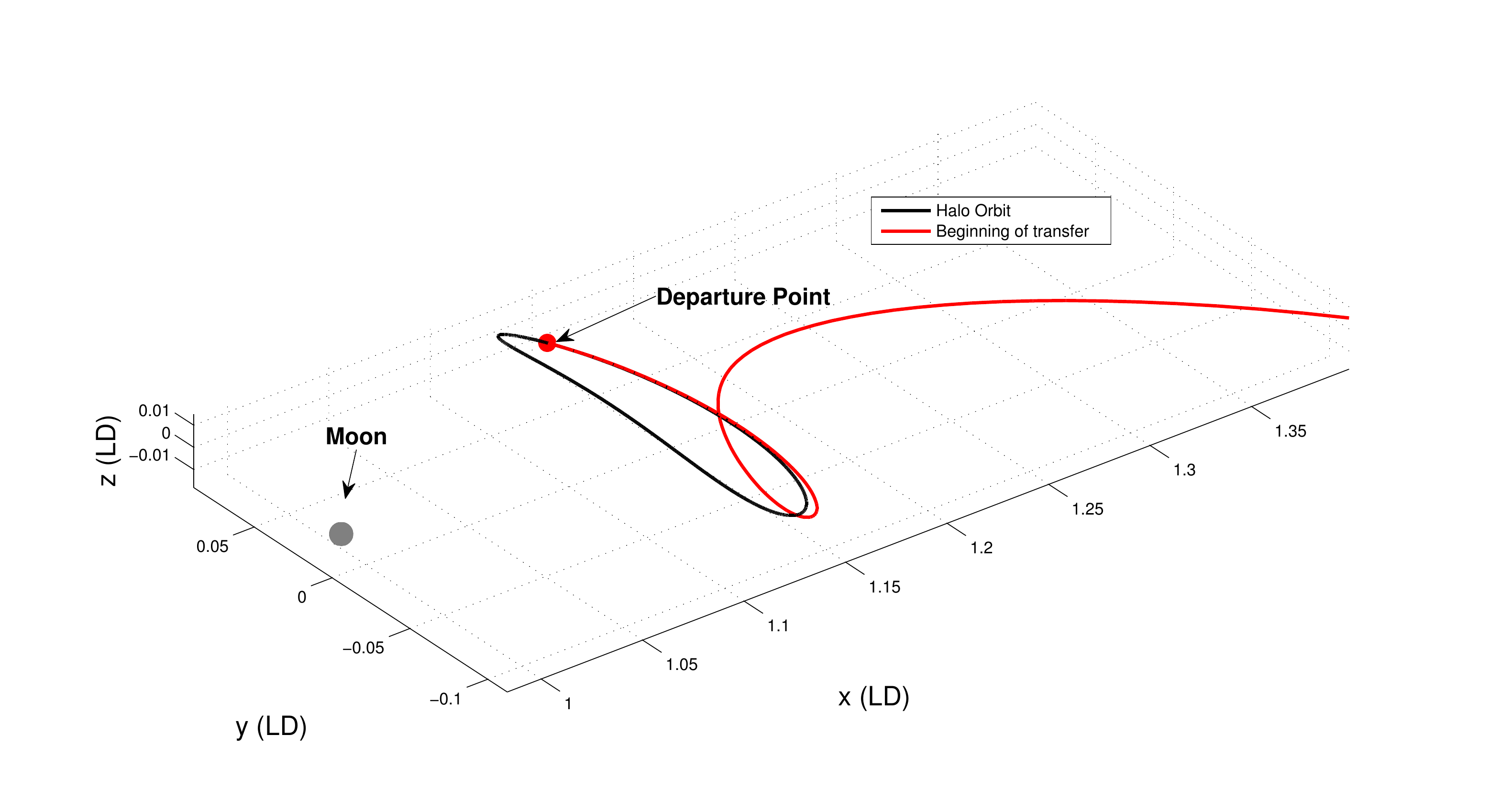}
\end{center}
\vspace*{-0.5cm}
\caption{\label{HaloL2_traj1} Best 3 boosts rendezvous transfer to \RH from a Halo orbit around $L_2$: Inertial Frame (left). Rotating Frame (Right). Bottom: zoom on the start of the departure from the Halo orbit (rotating frame).}
\end{figure}

\begin{figure}[ht!]
\centering
\includegraphics[width=6cm]{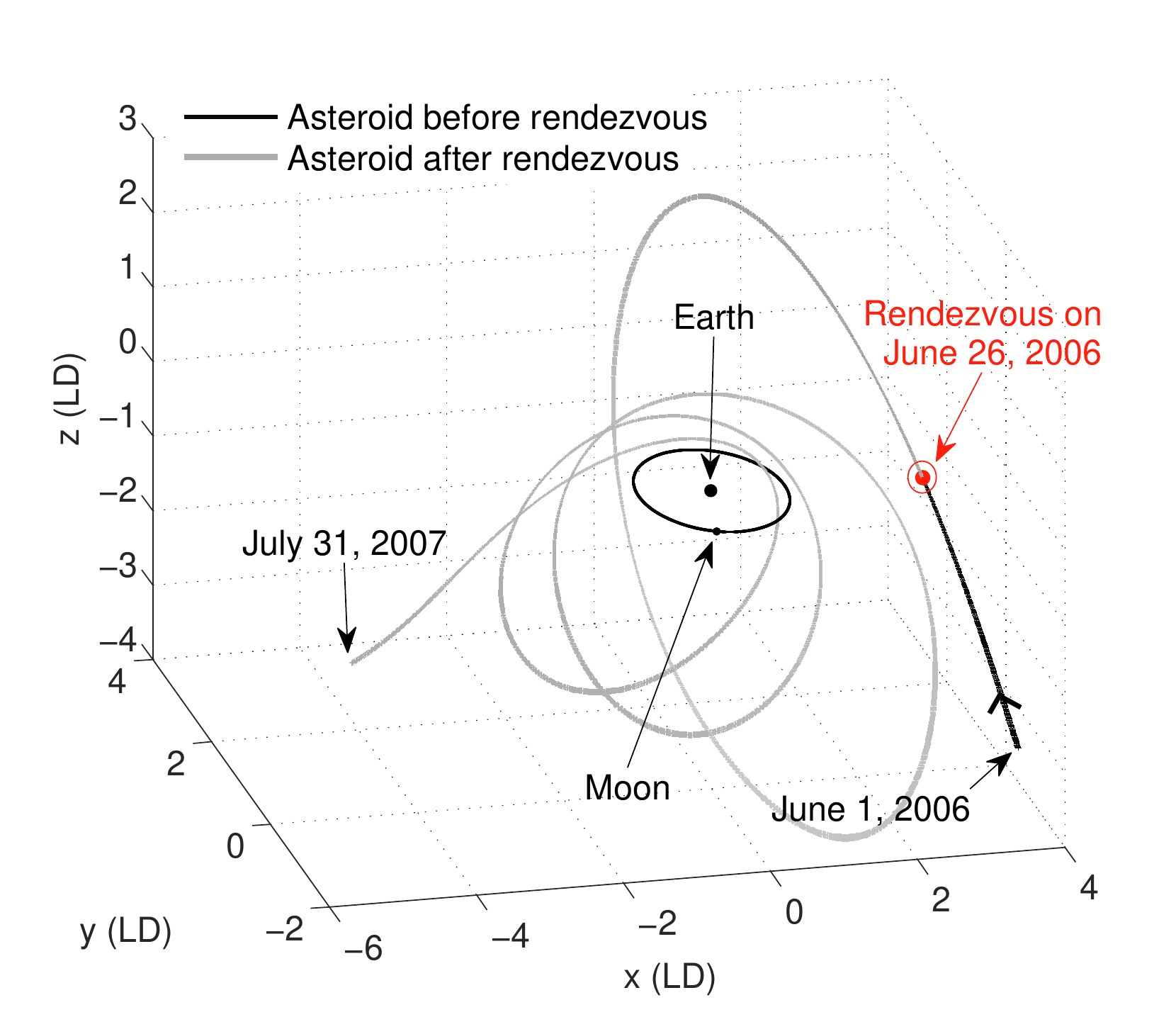}
\includegraphics[width=6cm]{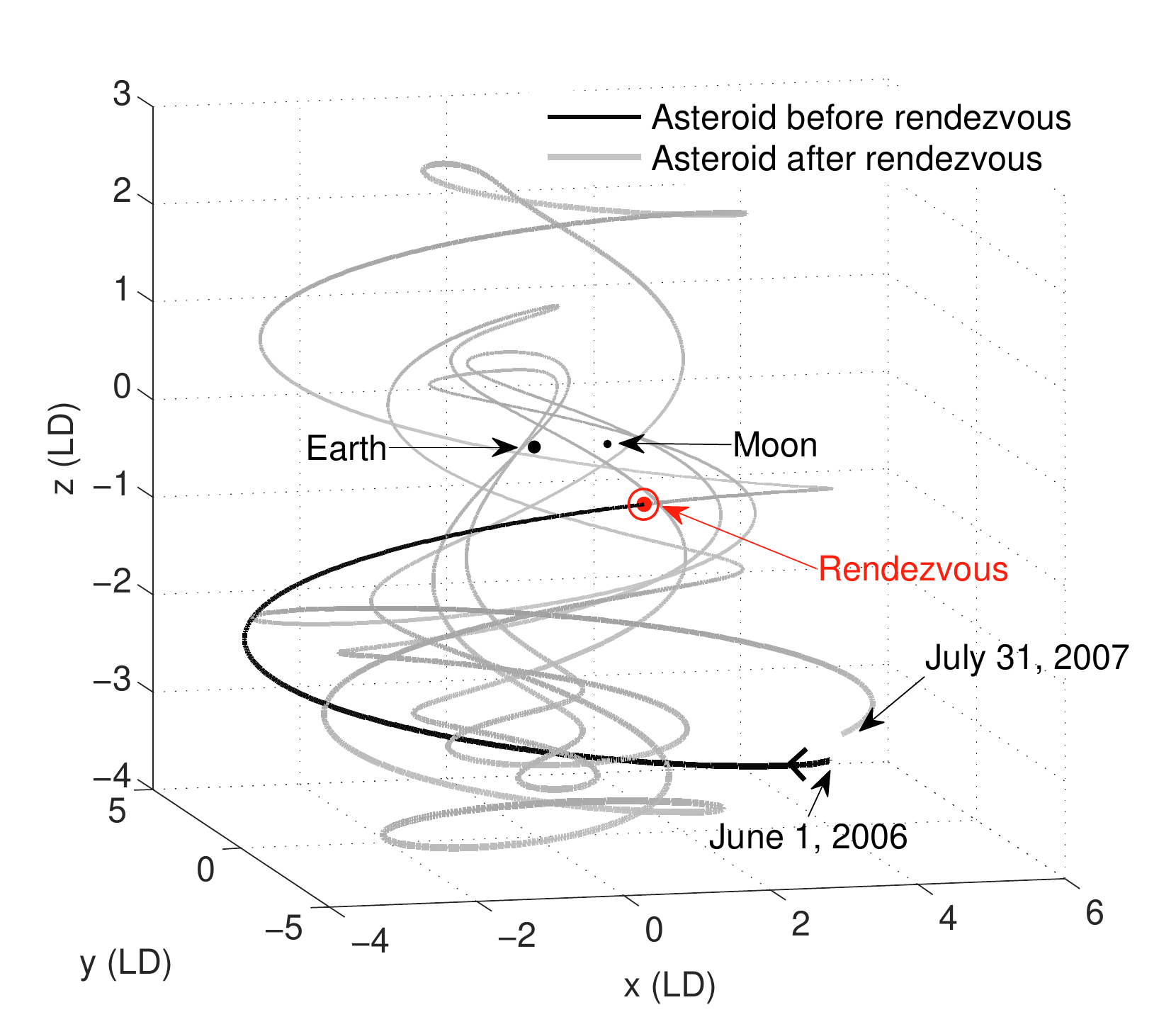}
\caption{Orbit of \RH in rotating (right) and inertial (left) reference frame with the rendezvous point corresponding to the best transfer. \label{RH120fulltraj-Halo}}
%\qquad
\end{figure}

\begin{figure}[ht!]
\begin{center}
\includegraphics[width=8cm]{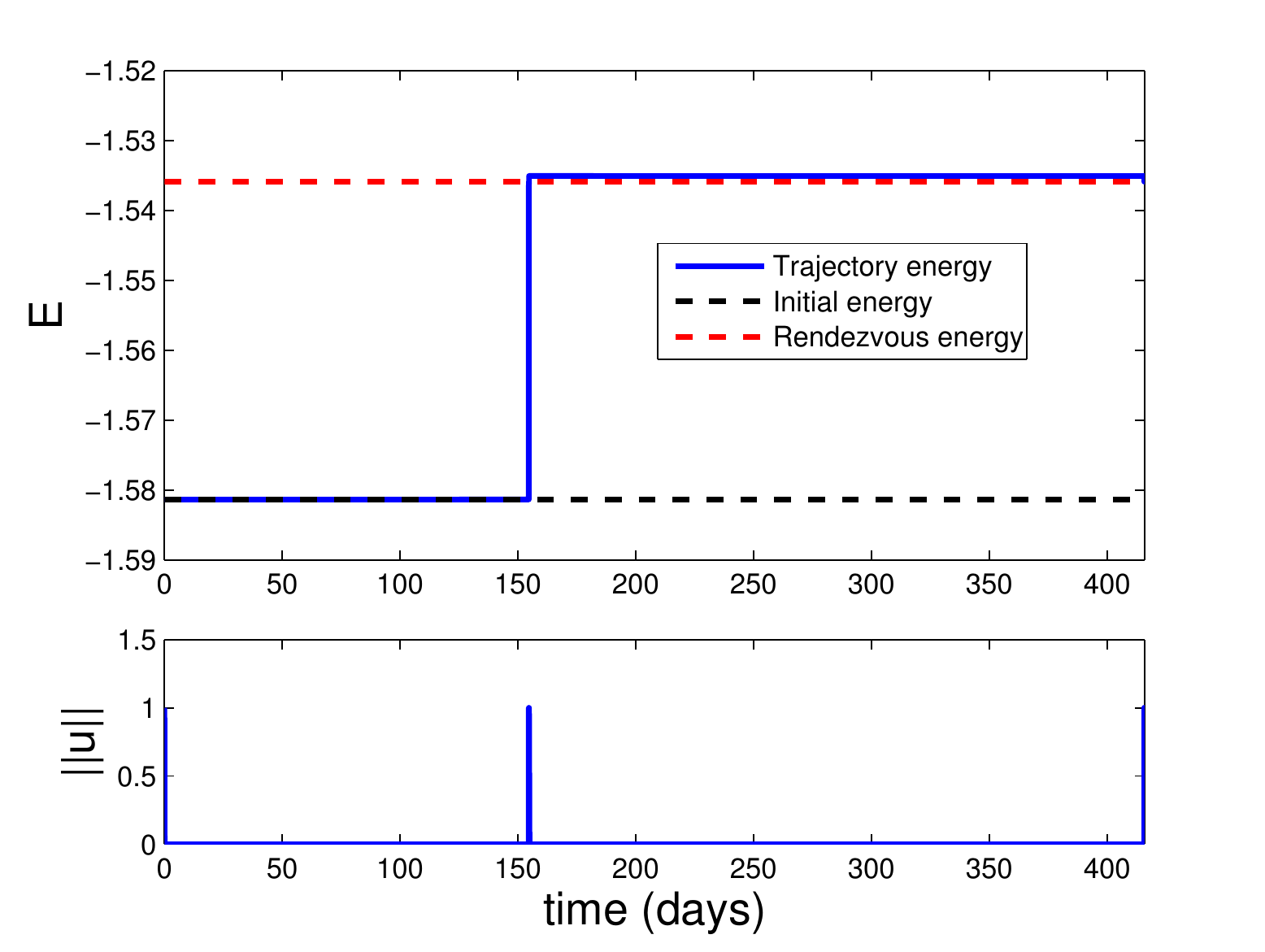}
\end{center}
\caption{\label{Control} The bottom graph represents the control strategy for the best transfer to \RH, and the top graph shows the evolution of the energy for this transfer. The energy of initial Halo orbit $E_{HaloL2}$ = -1.581360, the energy after 1st boost = -1.581335, which means that the first boost raised the energy for about $2.46e-05$. The energy after 2nd boost = -1.535022 which is an increase of about 0.046313 and the energy at the rendezvous point = -1.535804 which is a decrease of about = 7.814e-04. The total energy difference is about 0.0456.}
\end{figure}

The best transfer is represented on Figures \ref{HaloL2_traj1}, and some relevant data is also presented in Table \ref{Table-Best-Halo2-Transfer}. 
\begin{table}[ht!]
\begin{center}
\textbf{Table \ref{Table-Best-Halo2-Transfer}: Best Transfer Data from Halo-L$_\textnormal{2}$}
\begin{tabular}{@{}llll@{}}
\toprule
\textbf{Parameter} & \textbf{Symbol} & \textbf{Value} \\ \midrule
Transfer Duration& $t_f$ & 415.8 days \\
Final Mass & $m_f$ & 319.8 kg \\
Delta-v &  $\Delta v$ & 203.6 m/s \\
Final Position & $q_p^{\mathrm{rdv}}$ & \pbox{25cm}{(2.25, 3.21, -1.04)} \\
Final Velocities &  $q_v^{\mathrm{rdv}}$ & \pbox{25cm}{(2.92, -2.02, 0.46)}\\
Max Distance from L$ _2 $ & $d_{L_2}^\mathrm{max}$ & 10.63 LD \\ \bottomrule
\end{tabular}
\caption{\label{Table-Best-Halo2-Transfer} Data for the best transfer from $q_{\mathrm{HaloL}_2}$ to asteroid \RH.}
\end{center}
\end{table}
This transfer has a $\Delta v$ of $203.6\ m/s$. The rendezvous takes place on June 26$^{th}$ 2006 and lasts $415.5\ days$ which would require to detect and launch the mission about 14 months before June 1$^{st}$ 2006. 
On Figure \ref{RH120fulltraj-Halo}, we display the orbit of \RH in both the rotating and inertial frame with the rendezvous point for the best transfer. 
The best transfer exhibits 14 revolutions around the origin (in the rotating frame) and has a significant variation in the $z-$coordinate with respect to the EM plane. In particular, the $z-$coordinate of the rendezvous point is $-1.04$ normalized units, that is about 400 thousand $km$, and the maximum $z-$coordinate along the trajectory is $5.34$ normalized units, that is about 2 million $km$. The departure point on the parking orbit occurs $4.5$ days after $q_\mathrm{HaloL_2}$. The control strategy consists of a first boost of $19.7\ s$, a second boost starting after $154.7\ days$ and lasting $51.1\ min$ , finally, the last boost starting $261.1\ days$ after the second boost and lasts $13.8\ s$, see Figure \ref{Control}. 
Despite a seemingly short initial boost, this trajectory does not exhibit the profile of an initial jump on an unstable invariant manifold. However, this trajectory exploits the fact that a small initial boost leads to a far location where the gravity field of the two primaries is small and where the second boost can efficiently aim at the rendezvous point. In particular, we expect the existence of other local minima with a larger final time, going further away from the initial and final positions and providing an even better final mass. Also note that this kind of strategy could not have been obtained if we had restricted the control structure to have two boosts rather than three.

\begin{figure}[ht!]
\begin{center}
\includegraphics[width=7cm]{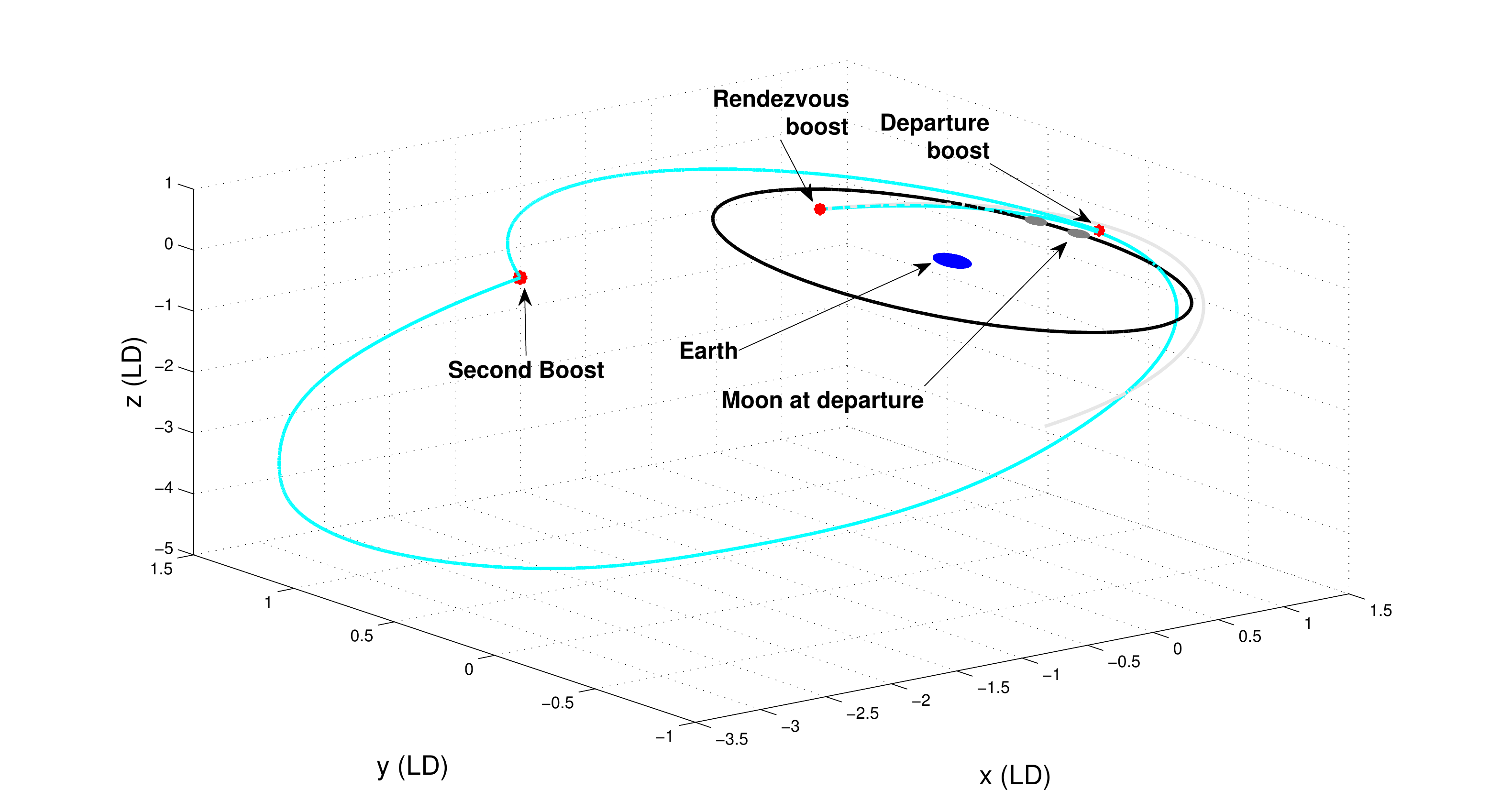}\includegraphics[width=7cm]{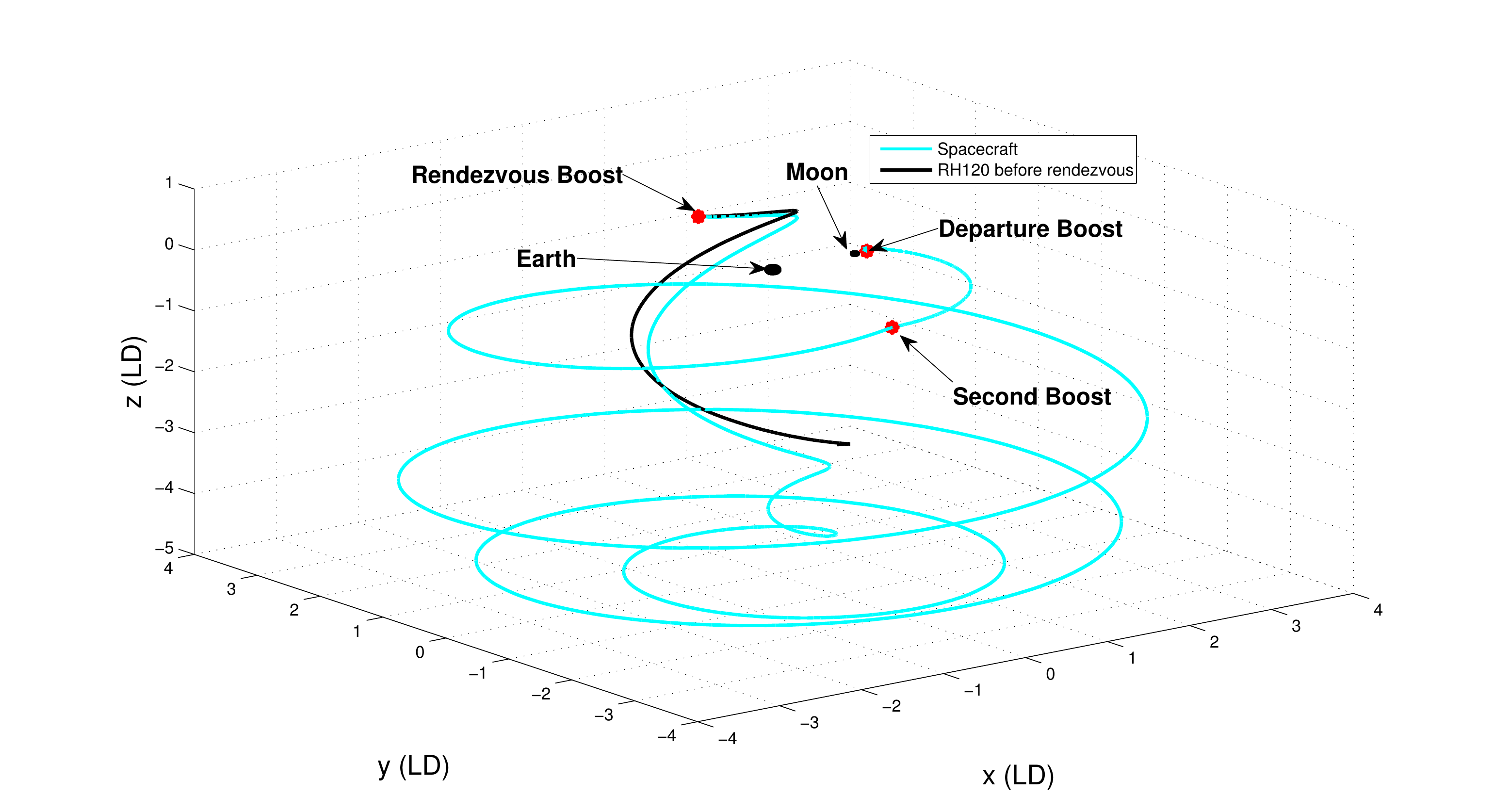}\\\includegraphics[width=9cm]{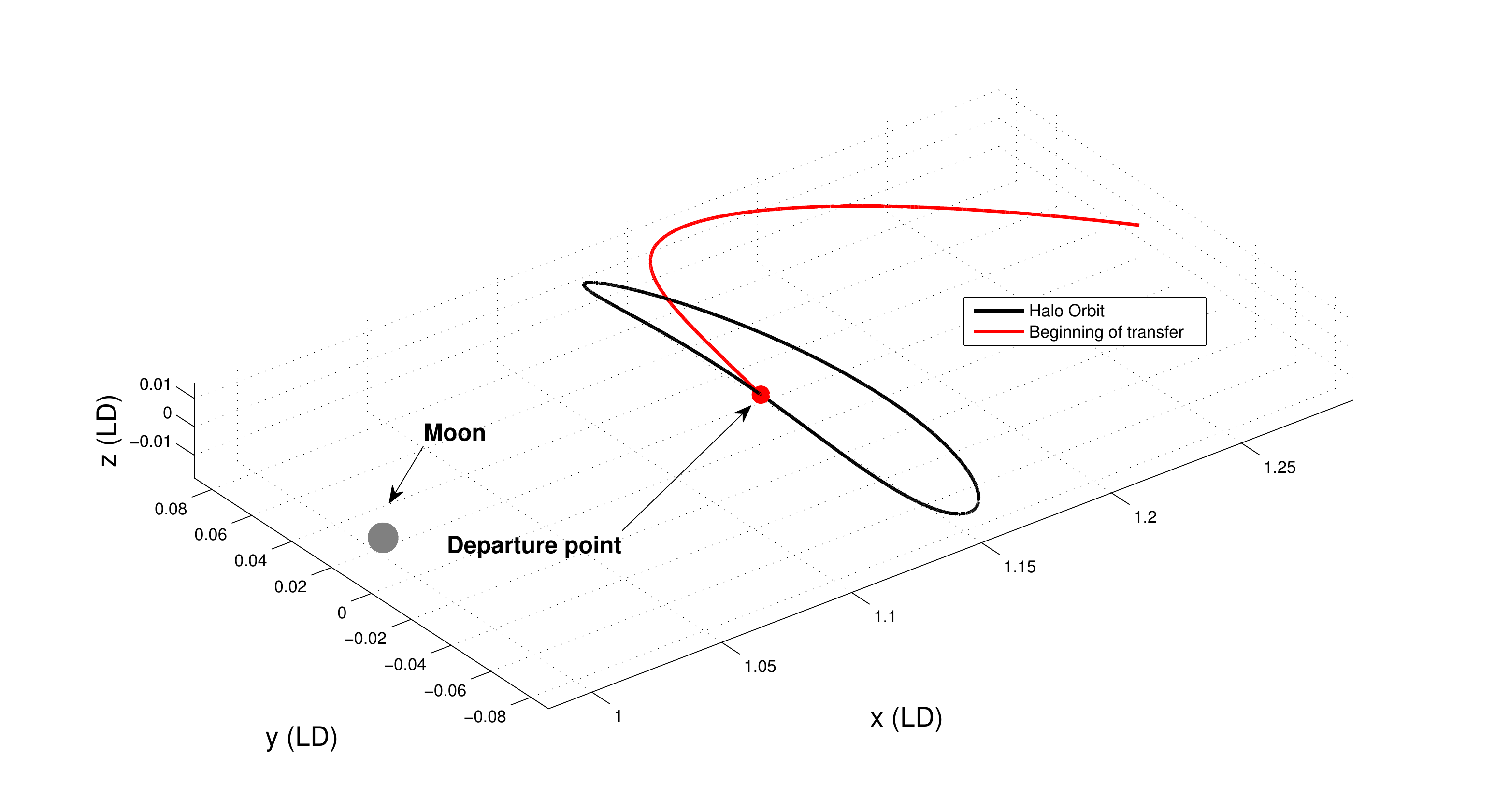}
\end{center}
\vspace*{-0.5cm}
\caption{\label{Worst} Worst 3 boosts rendezvous transfer to \RH from a Halo orbit around $L_2$: Inertial Frame (left). Rotating Frame (Right). Bottom: zoom on the start of the departure from the Halo orbit (rotating frame)}
\end{figure}

To contrast with the best transfer we represent the worst transfer in Figures \ref{Worst}. This transfer has a $\Delta v$ of $962.5\ m/s$ (final mass is $m_f=228.4 kg$) and a transfer time of $37.9$ days. The final position is $q_p^{\mathrm{rdv}}=(-0.65,0.35,0.92)$ and the final velocity is $q_v^{\mathrm{rdv}}=(-0.61,1.21,-0.28)$. The energy at the worst rendezvous point is about $-0.168$, and the maximum distance to $L_2$ is $6.21$ LD. 
\section{Conclusion and Future Work}
The numerical approach using an indirect method produces a low delta-v transfer to asteroid \RH departing from a Halo orbit around EM-$L_2$. Further calculations on four synthetic TCOs obtained from the databased produced in \cite{Jedicke} suggest that similar results can be expected on a large sample of TCOs. Indeed, in Table \ref{4TCO} we display data regarding the best three boosts transfer for four other TCOs. These transfers have a duration of about one year each and produce delta-v values between $223.9\ m/s$ and $344.2\ m/s$. The best transfer is for TCO$_1$ and occurs for an energy difference between terminal configurations of about $0.13$, for TCO$_{16}$ the difference is about $0.6$, $1.5$ for TCO$_{19}$ and $2.3$ for TCO$_{11}$. This reinforce the idea of a relation between the final mass of the transfer and the energy difference of the rendezvous point with respect to the one from the Halo orbit. This will be explored on a larger population of TCOs with the goal to characterize which asteroids are better suited for a rendezvous mission with a low delta-v. It can also be noted that the maximum distance of the spacecraft from $L_2$ during the transfer is similar for all four transfers which indicates that the long drift is used to pull away from the two primaries attraction fields to make the second boost more efficient. 
\begin{table}[ht!]
\begin{center}
\textbf{Table \ref{Table-Best-Halo2-Transfer-TCOs}: Best Transfer Data from Halo-L$_\textnormal{2}$ for Selected TCOs}
\begin{tabular}{@{}llllll@{}}
\toprule
\textbf{Parameter} & \textbf{Symbol} & \textbf{TCO$_1 $} & \textbf{TCO$_{11} $} & \textbf{TCO$_{16} $} & \textbf{TCO$_{19} $} \\ \midrule
Transfer Duration (days) & $t_f$ & 362.0 & 386.6 & 362.2 & 364.9 \\
Final Mass (kg) & $m_f$ & 316.9 & 300.5 & 311.0 & 307.1 \\
Delta-v (m/s) &  $\Delta v$ & 223.9 & 344.2 & 266.1 & 294.6 \\
Max Distance from L$ _2 $ & $d_{L_2}^\mathrm{max}$ & 12.7 & 11.5 & 11.5 & 12.8\\ 
Time to $d_{L_2}^\mathrm{max}$ (days)& $t(d_{L_2}^\mathrm{max})$ & 232.9 & 196.6 & 197.3 & 232.3 \\
\bottomrule
\end{tabular}
\caption{\label{Table-Best-Halo2-Transfer-TCOs} Best transfer data from $q_{\mathrm{HaloL}_2}$ to selected TCOs.}
\end{center}
\label{4TCO}
\end{table}

Forthcoming work will also add the Sun as a perturbation into the model, such as in \cite{Belbruno,Topputo}, to determine its impact on the transfers. Moreover, for practical reason it would be very interesting to consider the spacecraft parked on a Halo orbit around the ES $L_1$ libration point. 
\section{Acknowledgements}
We would like to specially thank Robert Jedicke and Mikael Granvik for their generous access to the database of synthetic TCOs, and in general, their help and support for this research. Geoff Patterson and Monique Chyba were partially supported by the National Science Foundation (NSF) Division of Mathematical Sciences, award \#1109937.





\end{document}